\documentclass[10pt]{amsart}

\usepackage{amssymb}
\usepackage{amscd}
\usepackage{upref}
\usepackage{tikz}
\theoremstyle{plain}
\newtheorem{thm}{Theorem}[section]
\newtheorem{cor}[thm]{Corollary}
\newtheorem{lem}[thm]{Lemma}
\newtheorem{prop}[thm]{Proposition}

\newtheorem*{Theorem A}{Theorem A}
\newtheorem*{Theorem B}{Theorem B}
\newtheorem*{Theorem C}{Theorem C}

\theoremstyle{definition}

\theoremstyle{remark}

\newtheorem{rem}[thm]{Remark}

\numberwithin{equation}{section}

\begin{document}
\title[TQFT String Operations in Open-Closed String Topology]
{TQFT String Operations in Open-Closed String Topology}
\author{Hirotaka Tamanoi}
\address[]
{Department of Mathematics, University of California Santa Cruz \newline
\indent Santa Cruz, CA 95064}
\email[]{tamanoi@math.ucsc.edu}
\date{}
\begin{abstract}
To open-closed cobordism surfaces, open-closed string topology associates topological quantum field theory (TQFT) operations, namely string operations, which depend only on homeomorphism types of surfaces and which satisfy the sewing property. We show that most TQFT string operations vanish in open-closed string topology. We describe those open-closed cobordisms with vanishing string operations, and give a short list of open-closed cobordisms with possibly nontrivial string operations. 
\end{abstract}
\maketitle

\tableofcontents

\section{Introduction}

In the open-closed string topology, (degree zero) string operations are associated to open-closed cobordism surfaces in such a way that string operations depend only on homeomorphism types of cobordism surfaces and they satisfy sewing property. Such open-closed string operations were constructed by Ramirez \cite{Ra} and Sullivan \cite{Su}. See also papers by Harrelsson \cite{Harrelsson} and by Lauda and Pfeiffer \cite{LP}. The open-closed string topology contains the closed string topology and string operations in closed string topology were constructed by Cohen and Godin \cite{CG}. In the same paper and in \cite{Go2}, they constructed higher string operations associated to homology classes of moduli spaces of Riemann surfaces with boundary, giving rise to a homological conformal field theory. The above mentioned string operations are associated to degree zero homology classes of the moduli spaces, and they give rise to topological quantum field theory (TQFT) \cite{At}. 

In \cite{T3}, we showed that higher genus degree zero closed string operations all vanish, and we described all cases of possibly nontrivial degree zero closed string operations. In this paper, we do the same for open-closed string topology.  As in the closed string topology case, most of the open-closed string operations vanish, and we will give a small and complete list of open-closed cobordisms whose associated string operations can be nontrivial. 

To describe the result, let $\Sigma$ be a connected open-closed cobordism, and suppose that end points of open strings are restricted to oriented closed submanifolds $I, J, K, L,\ldots$, so called D-branes. To $\Sigma$, we can associate the following numerical quantities. 
\begin{enumerate}
\item[] $g(\Sigma)=$ the genus of $\Sigma$. 
\item[] $\omega(\Sigma)=$ the number of windows (= free boundary circles) in $\Sigma$.  
\item[] $p(\Sigma)=$ the number of incoming closed strings. 
\item[] $q(\Sigma)=$ the number of outgoing closed strings. 
\item[] $r(\Sigma)=$ the number of boundary circles with only incoming open strings. 
\item[] $s(\Sigma)=$ the number of boundary circles with only outgoing open strings. 
\item[] $t(\Sigma)=$ the number of boundary circles containing both incoming and outgoing open strings. 
\end{enumerate} 

Let $P_{IJ}$ be the space of all continuous open strings in $M$ leaving $I$ and ending in $J$. Suppose $(I,J)$s are labels of end points of incoming open strings, and $(K,L)$s are labels of end points of outgoing open strings. For the construction of string operations, we require that there exists at least one outgoing closed string or an outgoing open string, namely, $q+s+t\ge1$. The associated string operation is of the form 
\begin{equation*}
\mu_{\Sigma}: H_*(LM)^{\otimes p}\otimes\bigotimes_{(I,J)} H_*(P_{IJ}) 
\longrightarrow H_*(LM)^{\otimes q}\otimes\bigotimes_{(K,L)} H_*(P_{KL}).
\end{equation*}
Our main theorem describes those open-closed cobordisms whose associated  degree zero string operations vanish. We assume that the oriented closed submanifolds $I,J,K,\dotsc,$ have dimension less that $d=\dim M$, although some statements are valid under weaker hypotheses. 

\begin{Theorem A} Let $\Sigma$ be a connected open-closed cobordism with at least one outgoing open or closed string, and let $\mu_{\Sigma}$ be the associated string operation. 

\noindent\textup{(I)} If $\Sigma$ satisfies one of the following conditions, then $\mu_{\Sigma}\equiv 0$. 
\begin{gather*}
\textup{(i)} \ \ g(\Sigma)\ge1. \qquad 
\textup{(ii)} \ \ \omega(\Sigma)\ge2. \qquad 
\textup{(iii)} \ \ t(\Sigma)\ge2. \\
\textup{(iv)} \ \ q(\Sigma)\ge 3. \qquad 
\textup{(v)} \ \ s(\Sigma)\ge1 \text{ and } s(\Sigma)+q(\Sigma)\ge2. 
\end{gather*}

\noindent\textup{(II)} Suppose $g(\Sigma)=0$ and $\omega(\Sigma)=1$. 
\begin{enumerate} 
\item[(i)] If $\Sigma$ has an incoming or outgoing open string $(r+s+t\ge1)$, then $\mu_{\Sigma}\equiv0$. 
\item[(ii)] If $\Sigma$ has no open strings $(r=s=t=0)$ and has at least two outgoing closed strings $(q\ge2)$, then $\mu_{\Sigma}\equiv0$. 
\end{enumerate} 

\noindent\textup{(III)} Suppose $g(\Sigma)=0$, $\omega(\Sigma)=0$ and $t(\Sigma)=1$. If $q+s\ge1$, then $\mu_{\Sigma}\equiv0$. 

\noindent\textup{(IV)} Suppose $g(\Sigma)=0$, $\omega(\Sigma)=0$, and $t(\Sigma)=0$. If $(q,s)\not=(1,0), (2,0), (0,1)$, then $\mu_{\Sigma}\equiv0$. 
\end{Theorem A} 

Part (IV) is included as a convenience, but it is a consequence of (iv) and (v) of part (I). To see this, since in string topology open-closed cobordism must have at least one outgoing open string or one outgoing closed string, when $t(\Sigma)=0$, we must have $q(\Sigma)+s(\Sigma)\ge1$. Then parts (iv) and (v) of Part (I) of Theorem A restrict the possible values of $(q,s)$ for nontrivial string operations to a set $(1,0)$, $(2,0)$, and $(0,1)$. This proves part (IV) of Theorem A.  

The proof of Theorem A is given in section 4, and it depends on the vanishing of string operations in the following four types proved in section 3. 
\begin{enumerate}
\item Vanishing of open window operation (Proposition \ref{open window operation}).
\item Vanishing of genus one operation (Proposition \ref{genus one operation}). 
\item Vanishing of handle attaching operation (Proposition \ref{handle attaching operation}). 
\item Vanishing of a double saddle operation (Proposition \ref{double saddle operation}). 
\item Vanishing of a double closed window operation (Proposition \ref{double closed window operation}).
\end{enumerate} 
The second and the third vanishing property above says that only genus $0$ open-closed cobordisms can have nontrivial string operations.
The first vanishing property above says that for a connected open-closed cobordism with windows, associated string operation vanishes unless all strings are closed, and in this case, there cannot be more than one window for non-triviality of the string operation by the fifth vanishing property. The fourth vanishing property implies that if a connected open-closed cobordism has at least two boundaries containing both incoming and outgoing open strings, then the associated string operation vanishes. 

Those open-closed cobordism surfaces $\Sigma$ not covered by Theorem A can have nontrivial string operations, and they can be classified into the following five different types. At the end of each types, we list their invariants.

\begin{Theorem B}
The following is a complete list of open-closed cobordisms with possibly nontrivial string operations.
\begin{enumerate} 
\item[(I)] $\Sigma_1$ is a closed cobordism of genus $0$ with exactly one window and with exactly one outgoing closed string. $(g=0,\omega=1, p\ge0, q=1, r=s=t=0)$. 
\item[(II)] $\Sigma_2$ is an open closed cobordism of genus zero with no windows with exactly one outgoing closed string and no outgoing open strings. $(g=\omega=0, p\ge0, q=1, r\ge0, s=t=0)$.
\item[(III)] $\Sigma_3$ is an open-closed cobordism of genus $0$ with no windows and no outgoing open strings, and has exactly two outgoing closed strings. $(g=\omega=0, p\ge0, q=2, r\ge0, s=t=0)$. 
\item[(IV)] $\Sigma_4$ is an open-closed cobordism of genus $0$ and no windows and no outgoing closed strings, and all outgoing open strings are along a single boundary of $\Sigma_4$. $(g=\omega=0, p\ge0, q=0, r\ge0, s=1, t=0)$. 
\item[(V)] $\Sigma_5$ is an open-closed cobordism of genus $0$ with no windows and no outgoing closed strings, and has exactly one boundary with both outgoing open strings and incoming open strings. $(g=\omega=0, p\ge0, q=0, r\ge0, s=0, t=1)$. 
\end{enumerate} 
\end{Theorem B}

Pictures of the above five types of surfaces are given below in figures 1 to 5 corresponding to (I) to (V) in Theorem B.

\begin{figure}
%Code for type (I) surface $Sigma_1$. 
\begin{center}
\begin{tikzpicture}
\path (0,0) coordinate (1);
%Multiple closed string product.
\draw (1) -- ++(-1,0) arc (270:225:1.5cm) arc (45:90:1.5cm) 
      -- ++(-0.5,0) node[] (4) {} arc (90:270:0.1cm and 0.27cm)
      -- ++(0.5,0) arc (90:0:0.5cm and 0.27cm) 
      arc (360:270:0.5cm and 0.27cm) -- ++(-0.5,0);
\draw (1)++(0,-1) -- ++(-1,0) arc (90:135:1.5cm) arc (315:270:1.5cm) 
      -- ++(-0.5,0) node[] (5) {} arc (270:90:0.1cm and 0.27cm) 
      -- ++(0.5,0) arc (270:360:0.5cm and 0.27cm) 
      arc (0:90:0.5cm and 0.27cm) -- ++(-0.5,0);
\draw[ultra thick] (4)++(0,-0.27) ellipse (0.1cm and 0.27cm);  
\draw[->,>=stealth, ultra thick] (4)++(0.1, -0.135) -- ++(0,0.1);   
\draw[ultra thick] (5)++(0,0.27)  ellipse (0.1cm and 0.27cm);
\draw[->,>=stealth, ultra thick] (5)++(0.1, 0.35) -- ++(0,0.1);
\draw[ultra thick] (1)++(-3.62,-0.5) ellipse (0.1cm and 0.27cm);
\draw[->,>=stealth,ultra thick] (1)++(-3.52,-0.37) -- ++(0,0.1);
%Closed window cobordism. 
\draw (1) arc (90:0:0.2cm and 0.5cm) arc (360:270:0.2cm and 0.5cm);
\draw[densely dashed] (1) arc (90:270:0.2cm and 0.5cm);
\draw (1)++(1,0) -- (1);
\draw (1)++(0,-1) -- ++(3.5,0);
\draw[ultra thick] (1)++(3.5,-1) arc (270:360:0.2cm and 0.5cm)
      arc (0:90:0.2cm and 0.5cm);
\draw[densely dashed,ultra thick] (1)++(3.5,0) arc (90:270:0.2cm and 0.5cm);
\draw (1)++(2.5,0) -- ++(1,0);
\path (1)++(2.5,0) coordinate (13);
\draw (1)++(1,0) to [out=270,in=270] (13);
\draw (1)++(1,0) to [out=330,in=210] (13);
\draw[->,>=stealth,ultra thick] (1)++(3.7,-0.4) -- ++(0,0.1);
\draw[->,>=stealth,ultra thick] (1)++(0.2,-0.4) -- ++(0,0.1);
\draw (1)++(4,-0.5) node[right,text width=3.5cm]
     {Exactly one outgoing closed string and no open strings ($q=1,s=t=0$).};
\draw (1)++(1.2,0.3) node[above,text width=6cm]
     {This genus $0$ cobordism has exactly one window ($g=0,w=1$).};
\draw (1)++(-4,-0.5) node[left,text width=3cm] {$p\ge0$ incoming closed strings, and no incoming open strings ($r=0$).};
\draw(1)++(0,-2) node[below,text width=8cm]
     {\textsc{Figure 1.} (I) Open-closed cobordism $\Sigma_1$.};
\end{tikzpicture}
\end{center}

%Code for type (II) surface $Sigma_2$. 
\begin{center}
\begin{tikzpicture}
\path (0,0) coordinate (3);
\foreach \r in {2.4}
{
\draw (3) -- ++(-1,0) arc (270:225:\r cm) arc (45:90:\r cm) 
      -- ++(-1,0);
\path (3) ++(-1,0) arc (270:225:\r cm) arc (45:90:\r cm) 
      ++(-1,0) coordinate (1);
\draw (3) ++(0,-1) -- ++(-1,0) arc (90:135:\r cm) 
      arc (315:270:\r cm) -- ++(-1,0) -- ++(0,1);
\path (3) ++(0,-1) ++(-1,0) arc (90:135:\r cm) 
      arc (315:270:\r cm) ++(-1,1) coordinate (2);
}
\path (3)++ (-3,-0.5) coordinate (0);
\draw (1) ++(0,-1) -- ++(1,0) to [out=0,in=90] (0)
      (2) -- ++(1,0)  to [out=0,in=270] (0);
\draw (1) arc (90:0:0.2cm and 0.5cm) arc (360:270:0.2cm and 0.5cm);
\draw[->,>=stealth,ultra thick] (1)++(0.2,-0.4) -- ++(0,0.1); 
\draw[densely dashed] (1) arc (90:270:0.2cm and 0.5cm);
\draw (2) arc (90:0:1cm and 0.5cm) arc (360:270:1cm and 0.5cm);
\draw[ultra thick, densely dashed] (3) arc (90:270:0.2cm and 0.5cm);
\draw[ultra thick] (3) arc (90:0:0.2cm and 0.5cm) 
      arc (360:270:0.2cm and 0.5cm);
%Add multiple loop product of closed strings. 
\draw (1) -- ++(-1,0) arc (270:225:1.5cm) arc (45:90:1.5cm) 
      -- ++(-0.5,0) node[] (4) {} arc (90:270:0.1cm and 0.27cm)
      -- ++(0.5,0) arc (90:0:0.5cm and 0.27cm) 
      arc (360:270:0.5cm and 0.27cm) -- ++(-0.5,0);
\draw (1)++(0,-1) -- ++(-1,0) arc (90:135:1.5cm) arc (315:270:1.5cm) 
      -- ++(-0.5,0) node[] (5) {} arc (270:90:0.1cm and 0.27cm) 
      -- ++(0.5,0) arc (270:360:0.5cm and 0.27cm) arc (0:90:0.5cm       and 0.27cm) -- ++(-0.5,0);
\draw[ultra thick] (4)++(0,-0.27) ellipse (0.1cm and 0.27cm);  
\draw[->,>=stealth, ultra thick] (4)++(0.1, -0.135) -- ++(0,0.1);   
\draw[ultra thick] (5)++(0,0.27)  ellipse (0.1cm and 0.27cm);
\draw[->,>=stealth, ultra thick] (5)++(0.1, 0.35) -- ++(0,0.1);
\draw[ultra thick] (1)++(-3.62,-0.5) ellipse (0.1cm and 0.27cm);
\draw[->,>=stealth,ultra thick] (1)++(-3.52,-0.37) -- ++(0,0.1);
%Add multiple open string product. 
\draw (2) -- ++(-1,0) arc (270:225:1cm) arc (45:90:1cm) 
      -- ++(-0.5,0);
\path (2) ++(-1,0) arc (270:225:1cm) arc (45:90:1cm) 
      ++(-0.5,0) coordinate (6);
\draw (2) -- ++(0,-1) -- ++(-1,0) arc (90:135:1cm) 
      arc (315:270:1cm) -- ++(-0.5,0) -- ++(0,0.75) -- ++(0.4,0) 
      arc (270:360: 0.4cm and 0.35cm) 
      arc (0:90:0.4cm and 0.35cm) -- ++(-0.4,0) ;
\path (2) ++(0,-1) ++(-1,0) arc (90:135:1cm) 
      arc (315:270:1cm) ++(-0.5,0) coordinate (7);
\path (7) ++(0,0.75) ++(0.4,0) arc (270:360: 0.4cm and 0.35cm) 
      arc (0:90:0.4cm and 0.35cm) ++(-0.4,0) coordinate (8);
\draw[ultra thick] (7) -- ++(0,0.75);
\draw[ultra thick] (8) -- (6);
\draw[->,>=stealth,ultra thick] (2)++(0,-0.4) -- ++(0,0.1);
\draw[->,>=stealth,ultra thick] (7)++(0,0.45) -- ++(0,0.1);
\draw[->,>=stealth,ultra thick]  (8)++(0,0.45) -- ++(0,0.1);
\draw (2)++(-3,-0.5) node[left,text width=4cm]
      {Incoming open strings along $r\ge0$ boundary circles.};
\draw[->,>=stealth,ultra thick] (3)++(0.2,-0.4) -- ++(0,0.01);
\draw (3)++(0.5,-1.5) node[below, text width=5cm]
     {Exactly one outgoing closed string and no outgoing open strings ($q=1,s=t=0$).};
\draw (1)++(-4,-0.5) node[left,text width=2.5cm] {$p\ge0$ incoming closed strings.};
\draw (1)++(1.5,0.3) node[above,text width=5cm]
       {This cobordism of genus $0$ has no windows ($g=0,w=0$).};
\draw (2)++(1,-1.8) node[below,text width=9cm]
      {\textsc{Figure 2.} (II) Open-closed cobordism $\Sigma_2$.};
\end{tikzpicture}
\end{center}

%Code for type (III) surface $Sigma_3$. 
\begin{center}
\begin{tikzpicture} 
\path (0,0) coordinate (3);
%Open-closed product.
\foreach \r in {2.4}
{
\draw (3) -- ++(-1,0) arc (270:225:\r cm) arc (45:90:\r cm) 
      -- ++(-1,0);
\path (3) ++(-1,0) arc (270:225:\r cm) arc (45:90:\r cm) 
      ++(-1,0) coordinate (1);
\draw (3) ++(0,-1) -- ++(-1,0) arc (90:135:\r cm) 
      arc (315:270:\r cm) -- ++(-1,0) -- ++(0,1);
\path (3) ++(0,-1) ++(-1,0) arc (90:135:\r cm) 
      arc (315:270:\r cm) ++(-1,1) coordinate (2);
}
\path (3)++ (-3,-0.5) coordinate (0);
\draw (1) ++(0,-1) -- ++(1,0) to [out=0,in=90] (0)
      (2) -- ++(1,0)  to [out=0,in=270] (0);
\draw (1) arc (90:0:0.2cm and 0.5cm) arc (360:270:0.2cm and 0.5cm);
\draw[->,>=stealth,ultra thick] (1)++(0.2,-0.4) -- ++(0,0.1); 
\draw[densely dashed] (1) arc (90:270:0.2cm and 0.5cm);
\draw (2) arc (90:0:1cm and 0.5cm) arc (360:270:1cm and 0.5cm);
\draw[->,>=stealth,ultra thick] (3)++(0.2,-0.4) -- ++(0,0.1);
%Add multiple loop product of closed strings. 
\draw (1) -- ++(-1,0) arc (270:225:1.5cm) arc (45:90:1.5cm) 
      -- ++(-0.5,0) node[] (4) {} arc (90:270:0.1cm and 0.27cm)
      -- ++(0.5,0) arc (90:0:0.5cm and 0.27cm) 
      arc (360:270:0.5cm and 0.27cm) -- ++(-0.5,0);
\draw (1)++(0,-1) -- ++(-1,0) arc (90:135:1.5cm) arc (315:270:1.5cm) 
      -- ++(-0.5,0) node[] (5) {} arc (270:90:0.1cm and 0.27cm) 
      -- ++(0.5,0) arc (270:360:0.5cm and 0.27cm) arc (0:90:0.5cm       and 0.27cm) -- ++(-0.5,0);
\draw[ultra thick] (4)++(0,-0.27) ellipse (0.1cm and 0.27cm);  
\draw[->,>=stealth, ultra thick] (4)++(0.1, -0.135) -- ++(0,0.1);   
\draw[ultra thick] (5)++(0,0.27)  ellipse (0.1cm and 0.27cm);
\draw[->,>=stealth, ultra thick] (5)++(0.1, 0.35) -- ++(0,0.1);
\draw[ultra thick] (1)++(-3.62,-0.5) ellipse (0.1cm and 0.27cm);
\draw[->,>=stealth,ultra thick] (1)++(-3.52,-0.37) -- ++(0,0.1);
%Add multiple open string product. 
\draw (2) -- ++(-1,0) arc (270:225:1cm) arc (45:90:1cm) 
      -- ++(-0.5,0);
\path (2) ++(-1,0) arc (270:225:1cm) arc (45:90:1cm) 
      ++(-0.5,0) coordinate (6);
\draw (2) -- ++(0,-1) -- ++(-1,0) arc (90:135:1cm) 
      arc (315:270:1cm) -- ++(-0.5,0) -- ++(0,0.75) -- ++(0.4,0) 
      arc (270:360: 0.4cm and 0.35cm) 
      arc (0:90:0.4cm and 0.35cm) -- ++(-0.4,0) ;
\path (2) ++(0,-1) ++(-1,0) arc (90:135:1cm) 
      arc (315:270:1cm) ++(-0.5,0) coordinate (7);
\path (7) ++(0,0.75) ++(0.4,0) arc (270:360: 0.4cm and 0.35cm) 
      arc (0:90:0.4cm and 0.35cm) ++(-0.4,0) coordinate (8);
\draw[ultra thick] (7) -- ++(0,0.75);
\draw[ultra thick] (8) -- (6);
\draw[->,>=stealth,ultra thick] (2)++(0,-0.4) -- ++(0,0.1);
\draw[->,>=stealth,ultra thick] (7)++(0,0.45) -- ++(0,0.1);
\draw[->,>=stealth,ultra thick]  (8)++(0,0.45) -- ++(0,0.1);
%Add closed string coproduct.    
\draw  (3) -- ++(1,0) arc (270:315:1cm) arc (135:90:1cm) -- ++(0.5,0)
      arc (90:0:0.1cm and 0.36cm) arc (360:270:0.1cm and 0.36cm) --                 ++(-0.5,0) arc (90:270:0.36cm) -- ++(0.5,0);
\path (3) ++(1,0) arc (270:315:1cm) arc (135:90:1cm) ++(0.5,0) 
       coordinate (14);
\draw (3)++(0,-1) -- ++(1,0) arc (90:45:1cm) arc (225:270:1cm) -- ++(0.5,0);
\path (3)++(0,-1) ++(1,0) arc (90:45:1cm) arc (225:270:1cm) ++(0.5,0)
       coordinate (15);
\draw[ultra thick, densely dashed] 
     (14) arc (90:270:0.1cm and 0.36cm);
\draw[ultra thick] (14) arc (90:0:0.1cm and 0.36cm) arc (360:270:0.1cm and          0.36cm);
\draw[->,>=stealth, ultra thick] (14)++(0.1,-0.3) -- ++(0,0.1);
\draw[ultra thick,densely dashed] (15) arc (270:90:0.1cm and 0.36cm);
\draw[ultra thick] (15) arc (270:360:0.1cm and 0.36cm) arc (0:90:0.1cm and 0.36cm);
\draw[->,>=stealth,ultra thick] (15)++(0.1,0.4) -- ++(0,0.1);
\draw[densely dashed] (3) arc (90:270:0.2cm and 0.5cm);
\draw (3) arc (90:0:0.2cm and 0.5cm) arc (360:270:0.2cm and 0.5cm);
\draw[->,>=stealth] (3)++(0.2,-0.4) -- ++(0,0.1);
%Text explanations. 
\draw (3)++(1.3,-1.8) node[below, text width=5cm]
     {Exactly two outgoing closed strings and no outgoing open strings ($q=2,s=t=0$).};
\draw (4)++(2,0.3) node[above,text width=5cm] 
     {$p\ge0$ incoming closed strings.};    
\draw (7)++(2,-0.3) node[below,text width=6cm]
      {Incoming open strings along $r\ge0$ boundary circles.};
\draw (3)++(0.5,0.3) node[above,text width=5cm]
      {This cobordism of genus $0$ has no windows ($g=\omega=0$).};
\draw (2)++(3,-2.7) node[below,text width=8cm]
      {\textsc{Figure 3.} (III) The open-closed cobordism $\Sigma_3$.};
\end{tikzpicture}
\end{center}

\end{figure}

\begin{figure}
%Code for type (IV) surface $Sigma_4$. 
\begin{center}
\begin{tikzpicture}
\path (0,0) coordinate (3); %The base point is left top corner. 
\draw (3) -- ++(1,0) arc (270:315:1.5cm) arc (135:90:1.5cm) -- ++(0.5,0)
      -- ++(0,-0.55) -- ++(-0.5,0) arc (90:270:0.5cm and 0.27cm) 
      -- ++(0.5,0);
\path (3) ++(1,0) arc (270:315:1.5cm) arc (135:90:1.5cm) ++(0.5,0)
       coordinate (9);
\path (9) ++(0,-0.55) ++(-0.5,0) arc (90:270:0.5cm and 0.27cm) 
      ++(0.5,0) coordinate (11);
\draw (3) -- ++(0,-1) -- ++(1,0) arc (90:45:1.5cm) arc (225:270:1.5cm) 
      -- ++(0.5,0) -- ++(0,0.55) -- ++(-0.5,0) 
      arc (270:90:0.5cm and 0.27cm) -- ++(0.5,0);
\path (3) ++(0,-1) ++(1,0) arc (90:45:1.5cm) arc (225:270:1.5cm) 
      ++(0.5,0) coordinate (10); 
\path (10) ++(0,0.55) ++(-0.5,0) arc (270:90:0.5cm and 0.27cm) 
      ++(0.5,0) coordinate (12);
\draw (3) -- ++(0,-1);
\draw[ultra thick] (9) -- ++(0,-0.55);
\draw[ultra thick] (10) -- ++(0,0.55);
\draw[ultra thick] (11) -- (12);
\draw[->,>=stealth,ultra thick] (9)++(0,-0.17) -- ++(0,0.1);
\draw[->,>=stealth,ultra thick] (10)++(0,0.35) -- ++(0,0.1);
\draw[->,>=stealth,ultra thick] (12)++(0,0.35) -- ++(0,0.1);
%Mouth piece. 
\draw (3) arc (90:270:1cm and 0.5cm);
\draw (3) -- ++(-2,0) 
      ++(0,-1) -- ++(2,0);
\draw[->, >=stealth, ultra thick] (3)++(0,-0.5) -- ++(0,0.1);
\draw[ultra thick] (3)++(-2,-0.5) ellipse (0.2cm and 0.5cm);
\draw[->, >=stealth, ultra thick] (3)++(-1.8,-0.5) -- ++(0,0.1);
\draw (3)++(4,-0.5) node[right,text width=3cm]
     {All outgoing open strings are along the same boundary circle.};
\draw (3)++(1.7,0.6) node[above,text width=10cm]
     {This cobordism has exactly two boundary circles, one of which is an incoming closed string and the other contains only outgoing open strings ($q=0,s=1,t=0$).}; 
\draw (3)++(2.5,-2.1) node[below,text width=12cm]
     {\textsc{Figure 4.} (IV) The open closed cobordism $\Sigma_4$ is obtained by sewing $\Sigma_2$ with the above cobordism along the outgoing closed string of $\Sigma_2$ and the above incoming closed string.};
\end{tikzpicture}
\end{center}

%Code for type (V) surface $Sigma_5$. 
\begin{center}
\begin{tikzpicture}
\path (0,0) coordinate (3);
\foreach \r in {2.2}
{
\draw (3) -- ++(-1,0) arc (270:225:\r cm) arc (45:90:\r cm) 
      -- ++(-1,0) -- ++(0,-1) -- ++(1,0);
\path (3) ++(-1,0) arc (270:225:\r cm) arc (45:90:\r cm) 
      ++(-1,0) coordinate (1);
\draw (3) -- ++(0,-1) -- ++(-1,0) arc (90:135:\r cm) 
      arc (315:270:\r cm) -- ++(-1,0) -- ++(0,1) -- ++(1,0);
\path (3) ++(0,-1) ++(-1,0) arc (90:135:\r cm) 
      arc (315:270:\r cm) ++(-1,0) ++(0,1) coordinate (2);
}
\path (3)++ (-3,-0.5) coordinate (0);
\draw (1)++(1,-1) to [out=0,in=90] (0)
      (2)++(1,0)  to [out=0,in=270] (0);
%Multiple open string coproduct. 
\draw (3) -- ++(1,0) arc (270:315:1.5cm) arc (135:90:1.5cm) -- ++(0.5,0)
      -- ++(0,-0.55) -- ++(-0.5,0) arc (90:270:0.5cm and 0.27cm) 
      -- ++(0.5,0);
\path (3) ++(1,0) arc (270:315:1.5cm) arc (135:90:1.5cm) ++(0.5,0)
       coordinate (9);
\path (9) ++(0,-0.55) ++(-0.5,0) arc (90:270:0.5cm and 0.27cm) 
      ++(0.5,0) coordinate (11);
\draw (3) -- ++(0,-1) -- ++(1,0) arc (90:45:1.5cm) arc (225:270:1.5cm) 
      -- ++(0.5,0) -- ++(0,0.55) -- ++(-0.5,0) 
      arc (270:90:0.5cm and 0.27cm) -- ++(0.5,0);
\path (3) ++(0,-1) ++(1,0) arc (90:45:1.5cm) arc (225:270:1.5cm) 
      ++(0.5,0) coordinate (10); 
\path (10) ++(0,0.55) ++(-0.5,0) arc (270:90:0.5cm and 0.27cm) 
      ++(0.5,0) coordinate (12);
\draw[ultra thick] (9) -- ++(0,-0.55);
\draw[ultra thick] (10) -- ++(0,0.55);
\draw[ultra thick] (11) -- (12);
\draw[->,>=stealth,ultra thick] (3)++(0,-0.4) -- ++(0,0.1);
\draw[->,>=stealth,ultra thick] (9)++(0,-0.17) -- ++(0,0.1);
\draw[->,>=stealth,ultra thick] (10)++(0,0.35) -- ++(0,0.1);
\draw[->,>=stealth,ultra thick] (12)++(0,0.35) -- ++(0,0.1);
%Mouth piece. 
\draw (1) arc (90:270:1cm and 0.5cm);
\draw (1) -- ++(-2,0) 
      ++(0,-1) -- ++(2,0);
\draw[->, >=stealth, ultra thick] (1)++(0,-0.5) -- ++(0,0.1);
\draw[ultra thick] (1)++(-2,-0.5) ellipse (0.2cm and 0.5cm);
\draw[->, >=stealth, ultra thick] (1)++(-1.8,-0.5) -- ++(0,0.1);
%Multiple open string product.
\draw (2) arc (270:225:1cm) arc (45:90:1cm) 
      -- ++(-0.5,0);
\path (2) arc (270:225:1cm) arc (45:90:1cm) 
      ++(-0.5,0) coordinate (6);
\draw (2) -- ++(0,-1) arc (90:135:1cm) 
      arc (315:270:1cm) -- ++(-0.5,0) -- ++(0,0.75) -- ++(0.4,0) 
      arc (270:360: 0.4cm and 0.35cm) 
      arc (0:90:0.4cm and 0.35cm) -- ++(-0.4,0) ;
\path (2) ++(0,-1) arc (90:135:1cm) 
      arc (315:270:1cm) ++(-0.5,0) coordinate (7);
\path (7) ++(0,0.75) ++(0.4,0) arc (270:360: 0.4cm and 0.35cm) 
      arc (0:90:0.4cm and 0.35cm) ++(-0.4,0) coordinate (8);
\draw[ultra thick] (7) -- ++(0,0.75);
\draw[ultra thick] (8) -- (6);
\draw[->,>=stealth,ultra thick] (2)++(0,-0.4) -- ++(0,0.1);
\draw[->,>=stealth,ultra thick] (7)++(0,0.45) -- ++(0,0.1);
\draw[->,>=stealth,ultra thick]  (8)++(0,0.45) -- ++(0,0.1);
\draw (1)++(3.5,0.3) node[above,text width=10cm]
     {This cobordism has exactly two boundary circles, one of which is a closed string to be sewn with $\Sigma_2$, and the other contains both incoming and outgoing open strings ($q=0,s=0,t=1$).};
\draw (2)++(-1.7,-1.7) node[below,text width=2cm]
     {Incoming open strings.};
\draw (3)++(3.5,1.6) node[below,text width=5cm]
     {Outgoing open strings.};
\draw (3)++(-0.5,-2.8) node[below,text width=10cm]
     {\textsc{Figure 5.} (V)  An open closed cobordism $\Sigma_5$ is obtained by sewing $\Sigma_2$ and the above cobordism along the closed string.}; 
\end{tikzpicture}
\end{center}

\end{figure}

String operations associated to type (I) surface $\Sigma_1$ is given purely in terms of loop products in $H_*(LM)$. The string operation associated to $\Sigma_2$ is given by open string products and closed string products, together with open-to-closed operation. String operations of type (III) and (IV) are computed from (II) by applying closed coproducts and open coproducts, together with closed-to-open operation. String operations of type (V) can be computed using the $H_*(LM)$-module structure on the path space homology $H_*(P_{IJ})$, and string operations associated to open-closed cobordisms homeomorphic to discs. These last disc operations will be examined in detail in \cite{T6}.

\bigskip

\section{Properties of transfer maps in homology and cohomology} 
%Section 2

We describe and prove various properties of transfer maps in the
context of fibration which will be useful in later sections and the subsequent paper \cite{T6}. 

First, we quickly review the construction of transfer maps. See \cite{CJ}. Let $M$ be a smooth closed oriented connected $d$-manifold and let $p:
E \rightarrow M$ be a Hurewicz fibration. For example, $E$ can be a
free loop space $LM$ of continuous maps from $S^1$ to $M$ and $p$ is an evaluation map at base points of loops, or a path space $P_{JK}$ consisting of continuous paths from
a submanifold $J$ to another submanifold $K$ and $p$ is an evaluation at time $0<t<1$. Let $\iota:L \rightarrow M$ be a smooth embedding of a closed oriented smooth
$\ell$-dimensional manifold, and let $q: E_L \rightarrow L$ be the
induced fibration. We use the same notation $\iota: E_L \rightarrow E$
to denote the inclusion map covering $\iota: L \rightarrow M$. The
following diagram provides the basic context.
\begin{equation*}
\begin{CD} 
E @<{\iota}<< E_L \\
@V{p}VV   @V{q}VV \\
M  @<{\iota}<< L 
\end{CD}
\end{equation*}
Let $\nu$ be the normal bundle to $\iota(L)$ in $M$, and we orient
$\nu$ so that we have an oriented isomorphism $\nu\oplus\iota_*(TL)\cong
TM|_{\iota(L)}$. Let $u\in \tilde{H}^{d-\ell}(L^{\nu})$ be the Thom
class of $\nu$. Let $N$ be a closed tubular neighborhood of $\iota(L)$
such that $D(\nu)\cong N$, where $D(\nu)$ is the closed disc bundle of
$\nu$ associated to some metric on $\nu$. We use the same notation for
the corresponding Thom class $u\in H^{d-\ell}(N,\partial N)$. Let
$\iota_N: N \rightarrow M$ be the inclusion map, and let $c: M
\rightarrow N/\partial N$ be the Thom collapse map. Let $v=c^*(u)\in
H^{d-\ell}(M)$ be the Thom class of the embedding $\iota$. The Thom
class $v$ is characterized by $v\cap [M]=\iota_*([L])$. In the
following commutative diagram, let $u'\in H^{d-\ell}(M,M-L)$ be the
corresponding Thom class, which is represented by a cocycle which
vanishes on singular simplices in $M$ not interesting with $L$.
\begin{equation*}
\begin{CD}
u'\in H^{d-\ell}(M, M-L) @>{j_M^*}>> H^{d-\ell}(M)\ni v  \\
@V{\iota_N^*}V{\cong}V  @A{c^*}AA  \\
u\in H^{d-\ell}(N, N-L) @>{j_N^*}>{\cong}> 
H^{d-\ell}(N,\partial N)\ni u
\end{CD}
\end{equation*}
We have $u=\iota_N^*(u')$ and $v=j_M^*(u')$. 

Let $\widetilde{N}=p^{-1}(N)$ and $\tilde{c}: E \rightarrow
\widetilde{N}/\partial\widetilde{N}$ be the Thom collapse map. Here,
$\partial\widetilde{N}$ is defined simply as $p^{-1}(\partial N)$. We
let $\iota_{\tilde{N}}: \widetilde{N} \rightarrow E$ be the inclusion
map.

The projection map $\pi: N \rightarrow L$ is a deformation
retraction. Using the homotopy lifting property of the Hurewicz
fibration $p: E \rightarrow M$, the map $\pi$ can be lifted to a
deformation retraction $\tilde{\pi}: \widetilde{N} \rightarrow E_L$ as
follows. Let $H: N\times I \rightarrow M$ be the homotopy for the deformation
retraction such that $H(x,0)=x$ and $H(x,1)=\pi(x)\in L$ for
$x\in N$. By homotopy lifting property, there exists a homotopy
$\widetilde{H} : \widetilde{N}\times I \rightarrow E$ extending the
identity map $\widetilde{H}(\tilde{x},0)=\tilde{x}$ and such that
$H\circ(p\times 1)=p\circ\widetilde{H}$. Then the map
$\tilde{\pi}(\tilde{x})=\widetilde{H}(\tilde{x},1)$ is a lift of
$\pi$, and is a deformation retraction.  We remark that although
$\tilde{\pi}$ is a deformation retraction, the map $\tilde{\pi}$ may
not be onto, and $\tilde{\pi}$ may not have a bundle structure. We can
construct a homotopy equivalence between the pull-back bundle
$q^*(\nu)$ and $\widetilde{N}$, but this map is in general far from
being a homeomorphism.

Let $\tilde{u}=p^*(u)\in H^{d-\ell}(\widetilde{N},
\partial\widetilde{N})$, $\tilde{v}=p^*(v)\in H^{d-\ell}(E)$, and
$\tilde{u}'=p^*(u')\in H^{d-\ell}(E, E-E_L)$ be pull-back Thom
classes. We consider the following commutative diagram, where
$\iota':E_L \rightarrow \widetilde{N}$ is an inclusion map and is a
homotopy inverse of $\tilde{\pi}: \widetilde{N} \rightarrow E_L$.
\begin{equation}\label{Thom classes}
\begin{CD}
\tilde{u}'\in H^{d-\ell}(E,E-E_L) @>{j_E^*}>> \tilde{v}\in H^{d-\ell}(E) @>{\iota^*}>> H^{d-\ell}(E_L) \\
@V{\iota_{\tilde{N}}^*}V{\cong}V @A{\tilde{c}^*}AA @A{{\iota'}^*}A{\cong}A \\
\tilde{u}\in H^{d-\ell}(\widetilde{N},\widetilde{N}-E_L) @>{j_{\tilde{N}}^*}>{\cong}> 
\tilde{u}\in H^{d-\ell}(\widetilde{N}, \partial\widetilde{N}) @>{j_{\tilde{N}}^*}>>
H^{d-\ell}(\widetilde{N})
\end{CD}
\end{equation}
Relations among various Thom classes are given by
$\tilde{v}=\tilde{c}^*(\tilde{u})$, $\tilde{v}=j_E^*(\tilde{u}')$,
$\tilde{u}=\iota_{\tilde{N}}^*(\tilde{u}')$, and
${\iota'}^*j_{\tilde{N}}^*(\tilde{u})=\iota^*(\tilde{v})$.

The homology and cohomology transfer maps $\iota_!$ and $\iota^!$ for
induced fibration are defined by the following compositions.
\begin{align*}
\iota_!&: H_*(E) \xrightarrow{\tilde{c}^*}
H_*(\widetilde{N},\partial\widetilde{N})
\xrightarrow{\tilde{u}'\cap(\ )} H_{*-d+\ell}(\widetilde{N})
\xrightarrow[\cong]{\tilde{\pi}_*} H_{*-d+\ell}(E_L), \\ \iota^!&:
H^*(E_L) \xrightarrow[\cong]{\tilde{\pi}^*} H^*(\widetilde{N})
\xrightarrow{\tilde{u}'\cup(\ )} H^{*+d-\ell}(\widetilde{N},
\partial\widetilde{N}) \xrightarrow{\tilde{c}^*} H^{*+d-\ell}(E).
\end{align*}
Here, the Thom maps given by cap or cup products with $\tilde{u}'$ may
not be an isomorphism since $\tilde{\pi}:\widetilde{N} \rightarrow
E_L$ does not have a bundle structure. However, for the rest of this
paper we do not need these Thom maps to be isomorphisms.

Note that if $L$ is not connected, then appropriate homology and
cohomology groups in the above diagrams split into direct sums.

By letting $p:E=M \rightarrow M$ be the trivial fibration, the following transfer maps between finite dimensional manifolds can be defined. 
\begin{equation*}
\iota_!:H_*(L) \longrightarrow H_{*-d+\ell}(L),\qquad
\iota^!:H^*(L) \longrightarrow H^{*+d-\ell}(M).
\end{equation*}

We collect basic properties of transfer maps which will be used later. See \cite{D} Chapter VIII \S10, \cite{Bredon} chapter VI \S14, and \cite{CoKl} for more on transfer maps. Although the properties discussed in the next theorem is well known, for convenience we provide their proofs, since orientation convention differs from literature to literature. 

\begin{thm}\label{transfer} Let $\iota:L \rightarrow M$ be a smooth 
embedding of a closed oriented $\ell$-manifold $L$ into a closed
oriented connected $d$-manifold $M$. Let $p:E \rightarrow M$ be a
Hurewicz fibration, and let $p:E_L \rightarrow L$ be the induced
fibration on $L$. Let $\iota: E_L \rightarrow E$ be the inclusion map.

The homology and cohomology transfer maps $\iota_!$ and $\iota^!$
satisfy the following identities for $a\in H_*(E_L)$, $b\in H_*(E)$,
$\alpha\in H^*(E)$, $\beta\in H^*(E_L)$, $\alpha'\in H^*(L)$.
\begin{align}
\iota_*\iota_!(b)&=\tilde{v}\cap b. 
\tag{1}\label{(1)}\\
\iota_!\iota_*(a)&=\iota^*(\tilde{v})\cap a=q^*(e_{\nu})\cap a. 
\tag{2}\label{(2)}\\
\iota_!(\alpha\cap b)
&=(-1)^{|\alpha|(d-\ell)} \iota^*(\alpha)\cap\iota_!(b). 
\tag{3}\label{(3)}\\
\iota^!\iota^*(\alpha)&=\tilde{v}\cup \alpha. 
\tag{4}\label{(4)}\\ 
\iota^*\iota^!(\beta)
&=\iota^*(\tilde{v})\cup\beta=q^*(e_{\nu})\cup\beta.   
\tag{5}\label{(5)}\\
\iota^!(\beta)\cap b
&=(-1)^{|\beta|(d-\ell)}\iota_*\bigl(\beta\cap\iota_!(b)\bigr).  
\tag{6}\label{(6)}\\ 
p^*\iota^!(\alpha')&=\iota^!\bigl(q^*(\alpha')\bigr).  
\tag{7}\label{(7)} \\
\iota_!([M])&=[L]. 
\tag{8}\label{(8)}\\
\iota^!(\{L_i\})&=(-1)^{\ell(d-\ell)}\{M\}, 
\tag{9}\label{(9)}
\end{align}
where in the last identity, $L=\coprod_iL_i$ is the decomposition into
connected components. Here $e_{\nu}=\iota^*(v)\in H^{d-\ell}(L)$ is
the Euler class of the normal bundle $\nu$ to $\iota(L)$ in $M$.
\end{thm} 
\begin{proof} Proofs are given in terms of commutative diagrams in 
terms of which the above identities become more or less transparent. 

(1) We consider the following commutative diagram.
\begin{equation*}
\begin{CD}
H_*(E) @>{\tilde{c}_*}>>  H_*(\widetilde{N},\partial\widetilde{N}) 
@>{\tilde{u}\cap(\ )}>>  H_{*-d+\ell}(\widetilde{N}) 
@>{\tilde{\pi}_*}>{\cong}> H_{*-d+\ell}(E_L) \\
@|  @V{{\iota_{\tilde{N}}}_*}V{\cong}V  @V{{\iota_{\tilde{N}}}_*}VV 
@V{\iota_*}VV   \\
H_*(E) @>{{j_E}_*}>>  H_*(E,E-E_L)  @>{\tilde{u}'\cap(\ )}>>  
H_{*-d+\ell}(E) @= H_{*-d+\ell}(E) 
\end{CD}
\end{equation*}
The middle square commutes because
$\tilde{u}=\iota_{\tilde{N}}^*(\tilde{u}')$ by \eqref{Thom classes}. The
composition of maps on the top line gives $\iota_!$. The commutativity
of the diagram implies for $b\in H_*(E)$,
\begin{equation*}
\iota_*\iota_!(b)=\tilde{u}'\cap (j_E)_*(b)=j_E^*(\tilde{u}')
\cap b=\tilde{v}\cap b. 
\end{equation*}
This proves (1). 

For (2), consider the following commutative diagram. 
\begin{equation*}
\begin{CD}
H_*(E_L) @>{\iota_*}>> H_*(E)  @= H_*(E)  \\
@|  @A{\iota_{\tilde{N}}}AA  @V{\tilde{c}_*}VV  \\
H_*(E_L) @>{\iota_*'}>{\cong}> H_*(\widetilde{N})  @>{{j_{\tilde{N}}}_*}>> 
H_*(\widetilde{N},\partial\widetilde{N})   \\
@V{\iota^*(\tilde{v})\cap(\ )}VV  @V{j_{\tilde{N}}^*(\tilde{u})\cap(\ )}VV 
@V{\tilde{u}\cap(\ )}VV   \\
H_*(E_L) @>{\iota_*'}>{\cong}> H_{*-d+\ell}(\widetilde{N}) @= H_{*-d+\ell}(\widetilde{N}) 
\end{CD}
\end{equation*}
Here, the bottom left square commutes since
${\iota'}^*j_{\tilde{N}}^*(\tilde{u})=\iota^*(\tilde{v})$ by
\eqref{Thom classes}. For $a\in H_*(E_L)$, the element $\iota_!\iota_*(a)$ is given
by following the diagram along the perimeter from the top left corner
to the bottom left corner in clockwise direction, noting that
$\iota_*'=(\tilde{\pi}_*)^{-1}$.  The commutativity of the diagram
then implies $\iota_!\iota_*(a)=\iota^*(\tilde{v})\cap a$.

(3) We consider the following commutative diagram, where $\alpha\in H^*(E)$. 
\begin{equation*}
\begin{CD}
H_*(E) @>{\tilde{c}_*}>> H_*(\widetilde{N},\partial\widetilde{N}) 
@>{\tilde{u}\cap(\ )}>> H_{*-d+\ell}(\widetilde{N}) @<{\iota_*'}<{\cong}< 
H_{*-d+\ell}(E_L)  \\
@V{\alpha\cap(\ )}VV  @V{\iota_{\tilde{N}}^*(\alpha)\cap(\ )}VV 
@V{\iota_{\tilde{N}}^*(\alpha)\cap(\ )}VV  @V{\iota^*(\alpha)\cap(\ )}VV \\
H_{*-|\alpha|}(E)  @>{\tilde{c}_*}>> 
H_{*-|\alpha|}(\widetilde{N},\partial\widetilde{N}) @>{\tilde{u}\cap(\ )}>> 
H_{*-|\alpha|-d+\ell}(\widetilde{N}) @<{\iota_*'}<{\cong}< H_{*-|\alpha|-d+\ell}(E_L)
\end{CD}
\end{equation*}
where the first square commutes since
$\tilde{c}_*=({\iota_{\tilde{N}}}_*)_*^{-1}\circ {j_E}_*$, and the
second square commutes up to $(-1)^{|\alpha|(d-\ell)}$ by an obvious
reason. The last square commutes because
$\iota=\iota_{\tilde{N}}\circ\iota'$. Since
$(\tilde{\pi}_*)^{-1}=\iota_*'$, the top and bottom rows both give
$\iota_!$.  Hence the commutativity of the diagram implies that
$\iota_!(\alpha\cap
b)=(-1)^{|\alpha|(d-\ell)}\iota^*(\alpha)\cap\iota_!(b)$.

(4) We consider the following commutative diagram.   
\begin{equation*}
\begin{CD}
H^*(E) @= H^*(E) @>{\tilde{u}'\cup(\ )}>> H^*(E,E-E_L)  @>{j_E^*}>> H^*(E) \\
@V{\iota^*}VV  @V{\iota_{\tilde{N}}^*}VV 
@V{\iota_{\tilde{N}}^*}V{\cong}V @|  \\
H^*(E_L) @>{\tilde{\pi}^*}>{\cong}>  H^*(\widetilde{N})
@>{\tilde{u}\cup(\ )}>>  H^{*+d-\ell}(\widetilde{N},\partial\widetilde{N}) 
@>{\tilde{c}^*}>> H^{*+d-\ell}(E)
\end{CD}
\end{equation*}
The first square commutes because $\iota\circ\tilde{\pi}\simeq
\iota_{\tilde{N}}$ since $\tilde{\pi}$ is a deformation
retraction. The second square commutes since
$\tilde{u}=\iota_{\tilde{N}}^*(\tilde{u}')$ by \eqref{Thom classes}. The
composition of maps in the bottom row is exactly $\iota^!$. Hence the
commutativity implies
$\iota^!\iota^*(\alpha)=j_E^*(\tilde{u}'\cup\alpha)
=j_E^*(\tilde{u}')\cup\alpha=\tilde{v}\cup\alpha$.

(5) We consider the following commutative diagram. 
\begin{equation*}
\begin{CD}
H^*(E_L) @<{{\iota'}^*}<{\cong}< H^*(\widetilde{N}) @>{\tilde{u}\cup(\ )}>> 
H^{*+d-\ell}(\widetilde{N},\partial\widetilde{N})  @>{\tilde{c}^*}>> 
H^{*+d-\ell}(E) \\
@VV{{\iota'}^*j_{\tilde{N}}^*(\tilde{u})\cup(\ )}V
@VV{j_{\tilde{N}}^*(\tilde{u})\cup(\ )}V  @VV{j_{\tilde{N}}^*}V  @| \\
H^{*+d-\ell}(E_L) @<{{\iota'}^*}<{\cong}< H^{*+d-\ell}(\widetilde{N}) 
@= H^{*+d-\ell}(\widetilde{N})   @<{\iota_{\tilde{N}}^*}<< H^{*+d-\ell}(E) 
\end{CD}
\end{equation*}
Here, the composition of maps in the top row gives $\iota^!$ and the
composition of maps in the bottom row gives $\iota^*$. Since
${\iota'}^*j_{\tilde{N}}^*(\tilde{u})=\iota^*(\tilde{v})$ by \eqref{Thom classes},
the commutativity of the diagram implies
$\iota^*\iota^!(\beta)={\iota'}^*j_{\tilde{N}}^*(\tilde{u})\cup \beta
=\iota^*(\tilde{v})\cup\beta$.

(6) Unraveling the definition of $\iota_!$, we have 
\begin{equation*}
\beta\cap\iota_!(b)
=(-1)^{|\beta|(d-\ell)}\tilde{\pi}_*\bigl((\tilde{u}\cup\tilde{\pi}^*(\beta))
\cap\tilde{c}_*(b)\bigr). 
\end{equation*}
Since $\iota\circ\tilde{\pi}\simeq\iota_{\tilde{N}}$ and 
$(\iota_{\tilde{N}})_*\circ\tilde{c}_*={j_E}_*$, we have 
\begin{align*}
(-1)^{|\beta|(d-\ell)}\iota_*\bigl(\beta\cap\iota_!(b)\bigr)
&={\iota_{\tilde{N}}}_*\bigl((\tilde{u}\cup\tilde{\pi}^*(\beta))
\cap\tilde{c}_*(b)\bigr) \\
&=(\iota_{\tilde{N}}^*)^{-1}\bigl(\tilde{u}\cup
\tilde{\pi}^*(\beta)\bigr)\cap(j_E)_*(b)  \\
&=(j_E)^*(\iota_{\tilde{N}}^*)^{-1}\bigl(\tilde{u}\cup
\tilde{\pi}^*(\beta)\bigr)\cap b.
\end{align*}
Since $(j_E)^*(\iota_{\tilde{N}}^*)^{-1}\bigl(\tilde{u}\cup
\tilde{\pi}^*(\beta)\bigr)=\tilde{c}^*\bigl(\tilde{u}\cup
\tilde{\pi}^*(\beta)\bigr)=\iota^!(\beta)$, the above formula 
is equal to $\iota^!(\beta)\cap b$. 

(7) This is straightforward. Since $\tilde{u}=p^*(u)$, we have 
\begin{equation*} \iota^!q^*(\alpha)=\tilde{c}^*
\bigl(\tilde{u}\cup\tilde{\pi}^*q^*(\alpha')\bigr)
=\tilde{c}^*\bigl(p^*(u)\cup p^*\pi^*(\alpha')\bigr)  
=p^*c^*(v\cup\pi^*(\alpha'))=p^*\iota^!(\alpha'). 
\end{equation*}

(8) By our choice of orientation on the normal bundle $\nu$, the Thom
isomorphism gives $\pi_*\bigl(u'\cap[N,\partial N]\bigr)=[L]$. Hence
we have
\begin{equation*}
\iota_!([M])=\pi_*\bigl(u'\cap c_*([M])\bigr)
=\pi_*(u'\cap[N,\partial N])
=[L].
\end{equation*}

(9) Let $L=\coprod_iL_i$ be the decomposition into path
components. Then for each component $L_i$, (6) and (8) implies
\begin{multline*}
\iota^!(\{L_i\})\cap[M]
=(-1)^{\ell(d-\ell)}\iota_*\bigl(\{L_i\}\cap\iota_!([M])\bigr) \\
=(-1)^{\ell(d-\ell)}\iota_*(\{L_i\}\cap[L])
=(-1)^{\ell(d-\ell)}[x_i],
\end{multline*}
where $x_i\in L_i$. Since $M$ is assumed to be connected,
$[x_i]=[x_0]$ for all $i$ for $x_0\in M$. Hence
$\iota^!(\{L\})\cap[M]=(-1)^{\ell(d-\ell)}|\pi_0(L)|[x_0]$. This is the same as $\iota^!(\{L\})=(-1)^{\ell(d-\ell)}|\pi_0(L)|\{M\}$.
\end{proof} 

\begin{rem}
If we change the orientation of the normal bundle $\nu$, then the sign
of the Thom class $u\in H^d(N,\partial N)$ also changes. Our choice of
the orientation is made so that we have
$\iota_!([M])=[L]$. Conversely, this identity characterizes the
orientation of $\nu$ and the sign of the Thom class $u$. 
\end{rem}

\bigskip

\section{Vanishing of some basic orientable open-closed string operations}
%Section 3

In this section, we first describe basic vanishing properties of open-closed string operations. Then in the next section, we give a refined list of open-closed cobordisms with vanishing string operations. In this process, we make a list of all those open-closed cobordisms with possibly nontrivial string operations. 

The proof of the vanishing of some basic string operations follows a similar pattern which we discuss first. Let $M$ be a connected oriented closed smooth $d$-manifold. Consider the following diagram where each square is a pull-back of Hurewicz fibrations, $\iota:L \rightarrow M$ is a smooth inclusion of an oriented closed submanifold $L$ of dimension $\ell$ which is not necessarily connected, and $\phi: L \rightarrow L\times L$ is the diagonal map: 
\begin{equation}\label{4-square pull-back}
\begin{CD}
P @<{\iota}<< Q @>{j}>> R @<{j}<< Q @>{\iota}>> P \\
@V{p}VV   @V{q}VV  @V{q'}VV  @V{q}VV  @V{p}VV \\
M @<{\iota}<< L @>{\phi}>> L\times L @<{\phi}<< L @>{\iota}>> M
\end{CD}
\end{equation}

\begin{lem}\label{basic lemma} In the homology diagram with transfers induced from the above diagram, for $a\in H_*(P)$, 
\begin{equation}
\iota_*j_!j_*\iota_!(a)=\chi(L)\bigl(p^*(\{M\})\cap a\bigr).
\end{equation}
\end{lem}
\begin{proof} Let $\nu$ be the normal bundle of the embedding $\phi$ oriented in such a way that we have an isomorphism of oriented vector bundles $\nu\oplus T\phi(L)\cong T(L\times L)|_{\phi(L)}$. Since $\nu\cong TL$ and its Euler class $e_{\nu}$ is given by $(-1)^{\ell}e_L=e_L$, using (2) of Theorem \ref{transfer}, we have 
$j_!j_*\iota_!(a)=q^*(e_L)\cap\iota_!(a)$. Now (6), (7) of Theorem \ref{transfer} implies 
\begin{align*}
\iota_*j_!j_*\iota_!(a)&=\iota_*\bigl(q^*(e_L)\cap \iota_!(a)\bigr) \\
&=(-1)^{\ell(d-\ell)}\iota^!\bigl(q^*(e_L)\bigr)\cap a \\
&=(-1)^{\ell(d-\ell)}p^*\bigl(\iota^!(e_L)\bigr)\cap a.
\end{align*}
Let $L=\coprod L_i$ be the decomposition into connected components of $L$. Then $e_L=\sum_ie_{L_i}=\sum_i\chi(L_i)\{L_i\}$. Using (9) of Theorem \ref{transfer}, we get 
\begin{equation*}
\iota^!(e_L)=\sum_i\chi(L_i)\iota^!(\{L_i\})
=\sum_I(-1)^{\ell(d-\ell)}\chi(L_i)\{M\}
=(-1)^{\ell(d-\ell)}\chi(L)\{M\}.
\end{equation*} 
Hence $\iota_*j_!j_*\iota_!(a)=\chi(L)\bigl(p^*(\{M\})\cap a\bigr)$. 
This completes the proof. 
\end{proof} 

We prove vanishing of five basic string operations using Lemma \ref{basic lemma}. Let $I,J,K,L$ be closed oriented submanifolds of $M$. Let $P_{IJ}$ be the space of continuous paths from points in $I$ to points in $J$. 

First, we consider a composition of an open string coproduct $\varphi_J$ and an open string product $\mu_J$: 
\begin{equation*}
\begin{CD}
H_*(P_{IK}) @>{\varphi_J}>> H_*(P_{IJ})\otimes H_*(P_{JK}) 
@>{\mu_J}>> H_*(P_{IK}).
\end{CD}
\end{equation*}
Since maps $\varphi_J$ and $\mu_J$ lower degree by $d-|J|$ and $|J|$, respectively, their composition $\mu_J\circ\varphi_J$ lowers degree by $d$. 
The relevant commutative diagram is the following one: 
\begin{equation*}
\begin{CD}
P_{IK} @<{\iota_J}<< P_{IJ}\underset{J}{\times}P_{JK} @>{j_J}>> P_{IJ}\times P_{JK} @<{j_J}<< P_{IJ}\underset{J}{\times}P_{JK} @>{\iota_J}>> P_{IK} \\
@V{p_{\frac12}}VV   @V{p}VV @V{p_1\times p_0}VV @V{p}VV @V{p_{\frac12}}VV   \\
M @<{\iota_J}<< J @>{\phi}>> J\times J @<{\phi}<< J @>{\iota_J}>> M 
\end{CD}
\end{equation*}
where $p_{\frac12}(\gamma)=\gamma(\frac12)$ for $\gamma\in P_{IK}$. The open string product and coproduct are given by $\mu_J=(\iota_J)_*(j_J)_!$ and $\varphi_J=(j_J)_*(\iota_J)_!$. We also use the notation $a\cdot_J b=(-1)^{j(|a|-j)}\mu_J(a\otimes b)$ for open string product for $a\in H_*(P_{IJ})$ and $b\in H_*(P_{JK})$. 

The corresponding open-closed cobordism $\Sigma$, which we call an open window, is an annulus with one incoming open string and one outgoing open string along the outer boundary, and its outer free boundaries carry labels $I,K$ and the inner open window carry a label $J$. See figure 6. 

%Open window operation. 
\begin{center}
\begin{tikzpicture}[>=stealth]
\draw (1,0) arc (0:90:0.5) -- ++(-1,0) arc (90:270:0.5) 
      -- ++(1,0) arc (270:360:0.5);
\draw (0,-1.5) -- ++(-0.5,0) arc (270:225:1.5) 
      arc (45:90:1.5) -- ++(-1,0) node[below] (1) {$I$};
\draw (0,-1.5) -- ++(0.5,0) arc (270:315:1.5) arc (135:90:1.5) 
      -- ++(1,0) node[below] (2) {$I$};
\draw (0,1.5) -- ++(-0.5,0) arc (90:135:1.5) 
      arc (315:270:1.5) -- ++(-1,0) node[above] (3) {$K$};
\draw (0,1.5) -- ++(0.5,0)  arc (90:45:1.5) arc (225:270:1.5) 
      -- ++(1,0) node[above] (4) {$K$};
\draw[ultra thick] 
(1) -- (3) (2) -- (4);
\draw[->, ultra thick] (1) -- (1)+(0,0.7);
\draw[->, ultra thick] (2) -- (2)+(0,0.7);
\path (0,-1.5) node[below] {$I$};
\path (0,1.5) node[above] {$K$};
\path (0,0.5) node[below] {$J$};
\path (0,-0.5) node[above] {$J$};
\draw[->] (0,0.5) -- (0,1);
\draw[->] (0,-1.5) -- (0,-1);
\draw (0,0.5)-- (0,1.5) (0,-1.5) -- (0,-0.5);
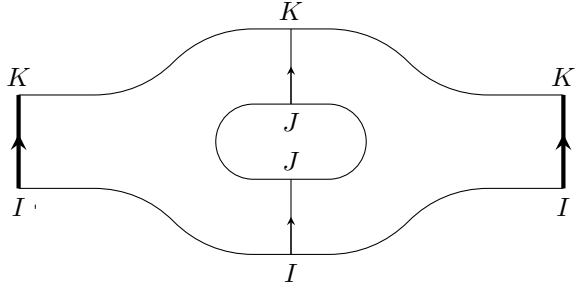
\node[text width=12cm, text centered] at (0,-2.5) {\textsc{Figure 6}. Open-window operation vanishes identically: 
$H_*(P_{IK}) \longrightarrow H_*(P_{IJ})\otimes H_*(P_{JK}) 
\longrightarrow H_*(P_{IK})$.};
\end{tikzpicture}
\end{center}

\begin{prop}[\textbf{Vanishing of open window operation}] \label{open window operation} Suppose $\dim I<\dim M$ or $\dim K<\dim M$. Then the open window operation $\mu_{\Sigma}$ vanishes. Namely, 
\begin{equation}
\mu_{\Sigma}=\mu_J\circ\varphi_J\equiv 0: H_*(P_{IK}) \longrightarrow H_*(P_{IJ})\otimes H_*(P_{JK}) 
\longrightarrow H_*(P_{IK}).
\end{equation}
\end{prop} 
\begin{proof} We give two different proofs. 

(Method I) For simplicity, we denote $\iota_J$ and $j_J$ simply by $\iota$ and $j$. By Lemma \ref{basic lemma}, for $a\in H_*(P_{IK})$, we have 
\begin{equation*}
\mu_J\circ\varphi_J(a)=\iota_*j_!j_*\iota_!(a)=
\chi(J)\bigl(p_{\frac12}^*(\{M\})\cap a\bigr).
\end{equation*}
Since $p_{\frac12}\simeq \iota_I\circ p_0: P_{IK} \rightarrow I \rightarrow M$, we have $p^*_{\frac12}(\{M\})=p^*_0\bigl(\iota_I^*(\{M\})\bigr)$. If $\dim I<\dim M$, then $\iota_I^*(\{M\})=0$, and consequently $p^*_{\frac12}(\{M\})=0$. Similarly, when $\dim K<\dim M$, we have $p^*_{\frac12}(\{M\})=p^*_1\iota_K^*(\{M\})=0$. Thus, in either case we have $\mu_J\circ\varphi_J(a)=0$ for every $a\in H_*(P_{IK})$. 

(Method II) Since $[I]\in H_*(P_{II})$ is the unit in the algebra $H_*(P_{II})$, we write $a\in H_*(P_{IK})$ as $a=[I]\cdot_I a$. Using the Frobenius property and associativity of open string products and coproducts, we have 
\begin{equation*}
\mu_J\circ\varphi_J(a)=\mu_J\circ\varphi_J([I]\cdot_I a) =\bigl(\mu_J\circ\varphi_J([I])\bigr)\cdot_I a.
\end{equation*}
Since the operator $\mu_J\circ\varphi_J$ lowers the degree by $d$, if $\dim I<\dim M$, then $\mu_J\circ\varphi_J([I])=0$. 
When $\dim K<\dim M$, by writing $a=a\cdot_K[K]$, we can argue as above to show that $\mu_J\circ\varphi_J(a)=0$. This completes the proof. 
\end{proof} 

Next we consider the closed string version of the above result. Let $\mu$ and $\varphi$ be the loop product and the loop coproduct: 
\begin{equation*}
\begin{CD}
H_*(LM) @>{\varphi}>> H_*(LM)\otimes H_*(LM) @>{\mu}>>
H_*(LM). 
\end{CD}
\end{equation*}
Both $\varphi$ and $\mu$ lower degree by $d$. The corresponding diagram is 
\begin{equation*}
\begin{CD}
LM @<{\iota}<<  LM\!\underset{M}{\times}\!LM @>{j}>>  LM\!\times \!LM
@<{j}<<  LM\!\underset{M}{\times}\!LM  @>{\iota}>> LM \\
@V{(p,p')}VV  @V{q}VV  @V{p\times p}VV   @V{q}VV  @V{(p,p')}VV  \\
M\!\times \!M   @<{\phi}<<  M  @>{\phi}>> M\!\times \!M  @<{\phi}<<  M 
@>{\phi}>> M\!\times \!M 
\end{CD}
\end{equation*}
where $p'(\gamma)=\gamma(\frac12)$. The loop product and coproduct are given by $\mu=\iota_* j_!$ and $\varphi=j_*\iota_!$. The open-closed cobordism $\Sigma$ associated to $\mu\circ\varphi$ is a torus with two boundary circles, one incoming and one outgoing closed strings. We call the associated string operation $\mu_{\Sigma}$ a genus $1$ operation. See figure 7. 

\begin{prop} [\textbf{Vanishing of genus one operation}] \label{genus one operation}  The string operation associated to a genus $1$ cobordism, a torus with one incoming and one outgoing closed strings, identically vanish. That is, 
\begin{equation}
\mu_{\Sigma}=\mu\circ\varphi\equiv 0: 
H_*(LM) \longrightarrow H_*(LM)\otimes H_*(LM) \longrightarrow H_*(LM).
\end{equation}
\end{prop} 
\begin{proof} For $a\in H_*(LM)$, by definition,  $\mu_{\Sigma}(a)=\iota_*j_!j_*\iota_!(a)$. Using Lemma \ref{basic lemma}, this is equal to 
\begin{equation*} 
\mu_{\Sigma}(a)=\chi(M)\bigl((p,p')^*(\{M\times M\})\cap a\bigr).
\end{equation*}
Here $(p,p')^*(\{M\times M\})=p^*(\{M\})\cup (p')^*(\{M\})=p^*(\{M\}\cup\{M\})=0$, since $p\simeq p'$ and $\{M\}\cup\{M\}=0$ in $H^*(M)$. Consequently, $\mu_{\Sigma}\equiv 0$. This completes the proof. 
\end{proof}

\begin{figure}
%Genus 1 operation.
\begin{center}
\begin{tikzpicture}[>=stealth]
\draw (0,0) ellipse (1.5cm and 1cm);
\draw (0,2) node[] (1) {} arc (90:150:3cm and 2cm) 
      arc (330:270:1cm and 0.66cm) --++(-0.5,0) node[] (2) {};
\draw (1) arc (90:30:3cm and 2cm) arc (210:270:1cm and 0.66cm) 
      -- ++(0.5,0) node[] (3) {};
\draw (0,-2) node[] (4) {} arc (270:210:3cm and 2cm) 
      arc (30:90:1cm and 0.66cm) -- ++(-0.5,0) node[] (5) {};
\draw (4) arc (270:330:3cm and 2cm) arc (150:90:1cm and 0.66cm) 
      -- ++(0.5,0) node[] (6) {};
\draw[densely dashed] (0,2) arc (90:270:0.2cm and 0.5cm);
\draw       (0,1) arc (270:360:0.2cm and 0.5cm) arc (0:90:0.2cm and 0.5cm);
\draw[->,>=stealth,ultra thick] (0.2,1.5) -- ++(0,0.1);
\draw[densely dashed] (0,-1) arc (90:270:0.2cm and 0.5cm);
\draw    (0,-2) arc (270:360:0.2cm and 0.5cm) arc (0:90:0.2cm and 0.5cm);
\draw[->,>=stealth,ultra thick] (0.2,-1.5) -- ++(0,0.1);
\draw[ultra thick] (-3.95,0) ellipse (0.25cm and 0.65cm);
\draw[->,ultra thick] (-3.7,0) -- (-3.7,0.1);
\draw[ultra thick] (3.95,0.65) arc (90:0:0.25cm and 0.65cm)
      arc (360:270:0.25cm and 0.65cm);
\draw[densely dashed,ultra thick] (3.95,0.65) 
      arc (90:270:0.25cm and 0.65cm);
\draw[->,ultra thick] (4.2,0) -- (4.2,0.1);
\path (0,-3) node[text width=10cm, text centered] 
{\textsc{Figure 7}. Genus $1$ operation vanishes identically:
$H_*(LM) \longrightarrow H_*(LM)\otimes H_*(LM) \longrightarrow H_*(LM)$.};
\end{tikzpicture}
\end{center}

\end{figure}

Note that the above genus one surface with two boundary circles can be obtained by attaching a handle to a cylinder. So we also consider an open-closed cobordism $\Sigma$ obtained by attaching a handle to a disc with one incoming and one outgoing open strings. Let the free boundary labels be submanifolds $I,J$. See figure 8. This open-closed cobordism describes a process in which a loop splits off from an open string and later rejoins to the open string. The following diagram describes this process. 
\begin{equation*}
\begin{CD}
P_{IJ} @<{\iota}<< P_{IJ}\!\underset{M}\times \!LM @>{j}>> P_{IJ}\!\times \!LM 
@<{j}<< P_{IJ}\!\underset{M}\times \!LM @>{\iota}>> P_{IJ} \\
@VV{(p_{\frac13},p_{\frac23})}V  @V{p}VV @V{p_{\frac12}\times p}VV 
@V{p}VV   @V{(p_{\frac13},p_{\frac23})}VV   \\
M\!\times \!M @<{\phi}<< M @>{\phi}>> M\!\times \!M @<{\phi}<< M @>{\phi}>> M\!\times\! M
\end{CD}
\end{equation*}
We call the associated string operation the handle attaching operation. 

\begin{prop}[\textbf{Vanishing of handle attaching operation}] \label{handle attaching operation} The string operation associated to attching a handle to a disc with one incoming and one outgoing open strings vanishes. That is, 
\begin{equation}
\mu_{\Sigma}\equiv 0: H_*(P_{IJ}) \longrightarrow H_*(P_{IJ})\otimes H_*(LM)
 \longrightarrow H_*(P_{IJ}).
\end{equation}
\end{prop} 
\begin{proof} For $a\in H_*(P_{IJ})$, using Lemma \ref{basic lemma}, we have 
\begin{equation*}
\mu_{\Sigma}(a)=\iota_*j_!j_*\iota_!(a)=\chi(M)\bigl((p_{\frac13},p_{\frac23})^*(\{M\times M\})\cap a\bigr).
\end{equation*}
Since the map $(p_{\frac13},p_{\frac23}):P_{IJ} \to M\times M$ is homotopic to a map $\phi\circ p_{\frac12}:P_{IJ} \to M\to M\times M$, we have $(p_{\frac13},p_{\frac23})^*(\{M\times M\})=p_{\frac12}^*\phi^*(\{M\times M\})=0$, since $\phi^*(\{M\times M\})=\{M\}\cup\{M\}=0\in H^*(M)$. This completes the proof. 
\end{proof}

\medskip

\begin{center}
%handle attaching to a disc
\begin{tikzpicture}[>=stealth]
%\draw[dashed] (2.5,0) arc (0:120:2.5cm and 1cm);
\draw (2.5,0) arc (0:60: 2.5cm and 1cm);
\draw (-2.5,0) arc (180:120:2.5cm and 1cm);
\draw (-2.5,0) arc (180:360:2.5cm and 1cm);
\draw[ultra thick] (2.5,0) arc (0:30:2.5cm and 1cm);
\draw[ultra thick] (-2.5,0) arc (180:210: 2.5cm and 1cm);
\draw[ultra thick] (-2.5,0) arc (180:150: 2.5cm and 1cm);
\draw[ultra thick] (2.5,0) arc (360:330:2.5cm and 1cm);
\draw (-1.5,0) arc (180:360:0.5cm and 0.2cm) 
arc (180:0:0.5cm) arc (180:360:0.5cm and 0.2cm) 
arc (0:180:1.5cm);
\draw (0,0.5) arc (270:360:0.15cm and 0.5cm) arc (0:90:0.15cm and 0.5cm);
\draw[dashed] (0,0.5) arc (270:90:0.15cm and 0.5cm);
\draw[->, ultra thick] (2.5,0) -- ++(0,0.1);
\draw[->, ultra thick] (-2.5,0) -- ++(0,0.1);
\path (1.8,-0.7) node[below] {$I$};
\path (1.8,0.7) node[above] {$J$};
\path (-1.8,0.7) node[above] {$J$};
\path (-1.8,-0.7) node[below] {$I$};
\draw[dashed] (1.5,0) arc(0:180:0.5cm and 0.2cm);
\draw[dashed] (-0.5,0) arc (0:180:0.5cm and 0.2cm);
\node[text width=12cm, text centered] at (0,-1.7) 
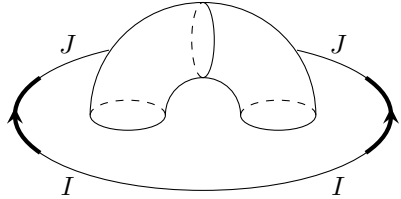
{\textsc{Figure 8}. Handle attaching operation vanishes identically:
$\mu\equiv 0:H_*(P_{IJ}) \longrightarrow H_*(P_{IJ})\otimes H_*(LM) \longrightarrow H_*(P_{IJ})$.};

\end{tikzpicture}
\end{center}

Propositions \ref{genus one operation} and \ref{handle attaching operation} show that all TQFT string operations associated to higher genus open-closed cobordism surfaces vanish. Thus only genus zero open-closed cobordisms can give rise to nontrivial TQFT string operations. Proposition \ref{genus one operation} was also proved in \cite{T3} in the context of closed string topology. The method used in that paper is slightly different from our current rather uniform method. In \cite{T3}, the coproduct $\varphi$ was determined explicitly, and $\varphi$ was shown to be nontrivial only on $H_d(LM)$ with values in $H_0(LM\times LM)$ generates by $[c_x]\otimes [c_x]$, where $[c_x]\in H_0(LM)$ is the homology class of the constant loop at $x\in M$. It was then shown that $\mu\circ\varphi\equiv 0$ since $\mu([c_x]\otimes[c_x])=0$ by dimensional reason.  The main result of this paper describes those open-closed cobordisms of genus 0 with vanishing string operations. Recall that in \cite{T3} we also showed that the string operation associated to a connected closed cobordism with at least three outgoing closed strings vanishes.  

Next, we discuss ``saddle'' interactions. This is the interaction of two open strings at their internal points. See \cite{T6} for more details. Let $I,J,K,L$ be closed oriented submanifolds of $M$. In the saddle interaction, an open string with end label $I,J$ and another open string with end label $K,L$ cut each other at their internal points and subsequently recombine. See figure 9.

\begin{figure}
%saddle interaction
\begin{center}
\begin{tikzpicture}[scale=0.7, auto] 
\path (-1,3) coordinate (1) node[above] {$I$};
\path (-3.42,4) coordinate (2) node[above] {$J$};
\path (-0.58,-4) coordinate (3) node[below] {$L$};
\path (-3,-3) coordinate (4) node[below] {$K$};
\path (4,3) coordinate (5) node[above] {$I$};
\path (1.58,4) coordinate (6) node[above] {$J$};
\path (4.42,-4) coordinate (7) node[below] {$L$};
\path (2,-3) coordinate (8) node[below] {$K$};
\path (-0.79,-3.25) coordinate (9);
\path (1.8,3.25) coordinate (10);
\path (2.53,0) coordinate (11);
\path (3.48,0) coordinate (12);
\draw (-2,-1) arc (270:360: 2 cm and 1 cm) 
arc (0:90:2 cm and 1 cm);
\draw[ultra thick] (1) parabola bend (-2,1) (2);
\draw[->, >=stealth, ultra thick] (-2,1) parabola  (-2.5,1.4);
\draw (2) to [out=350, in=200] node[above] {$J$} (6) ;
\draw (1) to [out=10, in=170] node[pos=0.45,above] {$I$} (5);
\draw[ultra thick] (4) parabola bend (-2,-1) (3);
\draw[->, >=stealth, ultra thick] (-2,-1) parabola (-1.5,-1.4);
\draw (3) to [out=10,in=170] node[below] {$L$}(7);
\draw[ultra thick] (7) to [out=110,in=250]  (5);
\path (2,-3) coordinate (8);
\draw (4) to [out=350,in=175] (9);
\draw[densely dashed] (9) to [out=355,in=190] node[pos=0.3,below] {$K$} (8);
\draw[densely dashed, ultra thick] (8) to [out=70,in=290] (10);
\draw[ultra thick] (10) -- (6);
\draw (0,0) arc (180:253:5 cm and 0.7 cm);
\draw[densely dashed] (0,0) arc (180:120:5 cm and 0.7 cm);
\draw[->, >=stealth, ultra thick] (11) -- (2.53,0.2);
\draw[->, >=stealth, ultra thick] (3.48,0.01) -- (12);
\draw node[text width=14cm, text centered] at (1,-6) 
{\textsc{Figure 9}. This open-closed cobordism is homeomorphic to a disc with $4$ open strings along the boundary, describing the saddle interaction of open strings:
$H_*(P_{IJ})\otimes H_*(P_{KL}) \longrightarrow H_*(P_{IL})\otimes H_*(P_{KJ})$.};
\end{tikzpicture}
\end{center}

\end{figure}

The saddle interaction is described by the following diagram: 
\begin{equation*}
\begin{CD}
P_{IJ}\times P_{KL} @<{\iota}<< P_{IJ}\underset{M}{\times}P_{KL} @>{j}>> P_{IL}\times P_{KJ}  \\
@VV{p_{\frac12}\times p_{\frac12}}V  @V{p}VV   @V{p_{\frac12}\times p_{\frac12}}VV \\
M\times M @<{\phi}<< M  @>{\phi}>> M\times M 
\end{CD}
\end{equation*}
Hence the saddle interaction or internal point interaction is given by the following composition of maps: 
\begin{equation*}
s=j_*\iota_!: H_*(P_{IJ})\otimes H_*(P_{KL}) \longrightarrow 
H_*(P_{IL})\otimes H_*(P_{KJ}).
\end{equation*}
This operation is associated to an open-closed saddle cobordism $\Sigma$ which is homeomorphic to a disc with four open strings such that incoming and outgoing open strings alternate along the boundary of the disc and its free boundaries are labeled by $I,J,K,L$ in this order. Saddle operations are not identically zero, but double saddle operations vanish identically.

\begin{prop} [\textbf{Vanishing of double saddle operation}] \label{double saddle operation} Suppose at least one of the closed oriented submanifolds $I,J,K,L$ have dimension less that $\dim M$. Then the double saddle interaction vanishes: 
\begin{equation} 
s\circ s\equiv 0: H_*(P_{IJ})\otimes H_*(P_{KL})  \rightarrow 
H_*(P_{IL})\otimes H_*(P_{KJ})  \rightarrow H_*(P_{IJ})\otimes H_*(P_{KL}).
\end{equation}
\end{prop}
\begin{proof} The diagram relevant to the double saddle interaction is the following. 
\begin{equation*}
\begin{CD}
P_{IJ}\!\times \!P_{KL} @<{\iota}<< P_{IJ}\!\underset{M}{\times}\!P_{KL} @>{j}>> P_{IL}\!\times \!P_{KJ}  @<{j}<<  P_{IL}\!\underset{M}{\times}\!P_{KJ}
@>{\iota}>> P_{IJ}\!\times \!P_{KL}\\
@V{p_{\frac12}\times p_{\frac12}}VV  @V{p}VV   @V{p_{\frac12}\times p_{\frac12}}VV  @V{p}VV   @V{p_{\frac12}\times p_{\frac12}}VV \\
M\!\times \!M @<{\phi}<< M  @>{\phi}>> M\!\times \!M  @<{\phi}<< M @>{\phi}>>  M\!\times \!M
\end{CD}
\end{equation*}
Here observe that intermediate spaces $P_{IJ}\underset{M}{\times}P_{KL}$ and $P_{IL}\underset{M}{\times} P_{KJ}$ are in fact exactly the same space, and we are in a position to apply Lemma \ref{basic lemma}. For $a\in H_*(P_{IJ})$ and $b\in H_*(P_{KL})$, we have 
\begin{align*}
(s\circ s)(a\times b)&=\iota_*j_!j_*\iota_!(a\times b)\\
&=\chi(M)\bigl[(p_{\frac12}\times 
p_{\frac12})^*(\{M\times M\})\cap(a\times b)\bigr] \\
&=(-1)^{|a|d}\chi(M)\bigl(p_{\frac12}^*(\{M\})\cap a\bigr)\times
\bigl(p_{\frac12}^*(\{M\})\cap b\bigr).
\end{align*}
If $\dim I<d$ or $\dim J< d$, then $p_{\frac12}^*(\{M\})=0\in H^*(P_{IJ})$, since $p_{\frac12}$ is homotopic to $\iota_I\circ p_0:P_{IJ} \to I \to M$, and also to $\iota_J\circ p_1:P_{IJ} \to J \to M$. Similarly, if $\dim K<d$ or $\dim L<d$, then $p_{\frac12}^*(\{M\})=0\in H^*(P_{KL})$. Hence in either case, we have $(s\circ s)(a\times b)=0$. This completes the proof. 
\end{proof} 

Next we consider the effect of a closed window labeled by an oriented closed submanifold $K$ of $\dim K=k$. Consider the following diagram: 
\begin{equation*}
\begin{CD}
LM @<{\iota}<< L_KM @>{j}>> P_{KK}  \\
@V{p}VV  @V{p}VV  @V{(p_0,p_1)}VV   \\
M @<{\iota_K}<<  K @>{\phi}>> K\times K
\end{CD}
\end{equation*}
Here $L_KM$ is the space of continuous loops in $M$ whose base points lie in $K$. Let $\theta_K=j_*\iota_!: H_*(LM) \rightarrow H_*(P_{KK})$ be an operation from closed strings to open strings lowering degree by $d-k$. It  describes an interaction in which a closed string touches the submanifold $K$ and splits into an open string whose end points lie in $K$. Similarly, the operation $\vartheta_K=\iota_*j_!: H_*(P_{KK}) \rightarrow H_*(LM)$ from open strings to closed strings lowering degree by $k$ describes an interaction in which open strings with end points in $K$ close up to become closed loops. The closed window operation $W_K$ is the composition of these two operations, and $W_K$ lowers degree by $d$:
\begin{equation*}
\begin{CD}
W_K=\vartheta_K\circ\theta_K: H_*(LM) \longrightarrow H_*(P_{KK}) \longrightarrow  H_*(LM).
\end{CD}
\end{equation*}
The open-closed cobordism $\Sigma$ for the string operation $W_K$ is a cylinder with one incoming and one outgoing closed strings and one hole, a free boundary circle, labeled by $K$. See figure 10. 

\medskip
%closed window operation
\begin{center}
\begin{tikzpicture}[>=stealth]
\draw[ultra thick] (-3,0) ellipse (0.3cm and 1cm);
\draw[ultra thick] (3,-1) arc (270:360:0.3cm and 1cm) 
arc (0:90:0.3cm and 1cm);
\draw[ultra thick, densely dashed] (3,1) arc (90:270:0.3cm and 1cm);
\draw[->,ultra thick] (-2.7,0) -- (-2.7,0.1);
\draw[->,ultra thick] (3.3,0) -- (3.3,0.1);
\draw (-3,-1) -- (3,-1) (-3,1) -- (-1.5,1) (1.5,1) -- (3,1);
\draw (-1.5,1) to [out=280, in=260] node[below] {$K$} (1.5,1) 
      (-1.5,1) to [out=330, in=210] node[above] {$K$} (1.5,1);
\draw (0,-1.5) node[below, text width=12cm, text centered] 
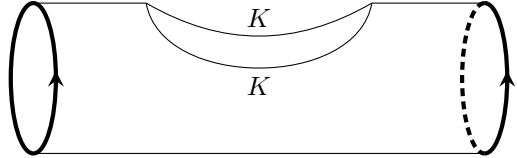
{\textsc{Figure 10}. Closed window operation of the form: 
$H_*(LM) \longrightarrow H_*(P_{KK}) \longrightarrow H_*(LM)$. };
\end{tikzpicture}
\end{center}

\begin{prop} [\textbf{Closed window operation}] For $a\in H_*(LM)$, the effect of a closed window operation $W_K$ with label $K$ is given by 
\begin{equation}
W_K(a)=\vartheta_K\circ\theta_K(a)=\chi(K)([c_0]\cdot a),
\end{equation}
where $[c_0]\in H_0(LM)$ is the homology class of a constant loop, and the dot $\cdot$ denotes the loop product in $H_*(LM)$. 
\end{prop}
\begin{proof} The diagram for the operation $W_K=\vartheta_K\circ\theta_K$ is the following one:  
\begin{equation*}
\begin{CD}
LM @<{\iota}<<  L_KM  @>{j}>>  P_{KK} @<{j}<< L_KM  @>{\iota}>>  LM \\
@V{p}VV  @V{p}VV   @V{(p_0,p_1)}VV  @V{p}VV  @V{p}VV   \\
M @<{\iota_K}<<  K  @>{\phi}>>  K\times K  @<{\phi}<<   K  @>{\iota_K}>>  M 
\end{CD}
\end{equation*}
For $a\in H_*(LM)$, we have $\vartheta_K\circ\theta_K(a)=\iota_*j_!j_*\iota_!(a)$. By Lemma \ref{basic lemma}, this is equal to 
\begin{equation*}
\vartheta_K\circ\theta_K(a)=\chi(K)\bigl(p^*(\{M\})\cap a\bigr)
=\chi(K)([c_0]\cdot a). 
\end{equation*}
Here we recall that for any cohomology class $\alpha\in H^*(M)$ and a homology class $b\in H_*(LM)$, we have $p^*(\alpha)\cap b=D(\alpha)\cdot b$, where $D(\alpha)\in H_*(M)$ is the Poincar\'e dual of $\alpha$. See  \cite{T2} for more details and related topics. This completes the proof. 
\end{proof} 

We can prove the same result using a different method. Let $\Sigma_1$ be an open-closed cobordism homeomorphic to a cylinder with one outgoing closed string and one free boundary circle labeled by $K$. The associated string operation $\mu_{\Sigma_1}: k \rightarrow H_*(LM)$, where $k$ is the ground field, is given by $\mu_{\Sigma_1}(1)=\vartheta_K([K])=\chi(K)[c_0]\in H_0(LM)$. Let $\Sigma'$ be an open closed cobordism obtained by sewing $\Sigma_1$ to one of the incoming closed string of a pair of pants representing the loop product. Then $\Sigma'$ is homeomorphic to the open-closed cobordism $\Sigma$ for the closed window operation $W_K$. See figure 11.

\begin{center}
\begin{tikzpicture}
\path (0,0) coordinate (3);
\foreach \r in {1.8}
{
\draw (3) -- ++(-1,0) arc (270:225:\r cm) arc (45:90:\r cm) 
      -- ++(-2,0) -- ++(0,-1); 
\path (3) ++(-1,0) arc (270:225:\r cm) arc (45:90:\r cm) 
      ++(-2,0) coordinate (1);
\draw (3) ++(0,-1) -- ++(-1,0) arc (90:135:\r cm) 
      arc (315:270:\r cm) -- ++(-2,0);
\path (3) ++(0,-1) ++(-1,0) arc (90:135:\r cm) 
      arc (315:270:\r cm) ++(-2,1) coordinate (2);
}
\draw[->,>=stealth] (1)++(0,-0.4) -- ++(0,0.01);
\path (3)++ (-3,-0.5) coordinate (0);
\draw (1) ++(0,-1) -- ++(2,0) to [out=0,in=90] (0)
      (2) -- ++(2,0)  to [out=0,in=270] (0);
\draw (2)++(1.5,0) arc (90:0:0.2cm and 0.5cm) arc (360:270:0.2cm and 0.5cm);
\draw[densely dashed] (2)++(1.5,0) arc (90:270:0.2cm and 0.5cm);
\draw[->,>=stealth] (2)++(1.7,-0.4) -- ++(0,0.01);
\draw (1)++(1.5,0) arc (90:0:0.2cm and 0.5cm) arc (360:270:0.2cm and 0.5cm);
\draw[densely dashed] (1)++(1.5,0) arc (90:270:0.2cm and 0.5cm);
\draw[->,>=stealth] (1)++(1.7,-0.4) -- ++(0,0.01);
\draw[ultra thick] (2) arc (90:360:0.2cm and 0.5cm) arc (0:90:0.2cm and 0.5cm); 
\draw[->,>=stealth, ultra thick] (2)++(0.2,-0.4) -- ++(0,0.01); 
\draw (1) arc (90:0:0.7cm and 0.5cm) arc (360:270:0.7cm and 0.5cm);
\draw (1) arc (90:270:0.7cm and 0.5cm);
\draw[ultra thick, densely dashed] (3) arc (90:270:0.2cm and 0.5cm);
\draw[ultra thick] (3) arc (90:0:0.2cm and 0.5cm) 
      arc (360:270:0.2cm and 0.5cm);
\draw[->,>=stealth,ultra thick] (3)++(0.2,-0.4) -- ++(0,0.01);

\draw (1)++(0.7,-0.5) node[left]{$K$};
\draw (1)++(-0.7,-0.5) node[left] {$K$};
\draw (1)++(1,0.1) node[above,text width=4cm]
       {This upper part is $\Sigma_1$.};
\draw (2)++(2.5,-1.3) node[below,text width=9cm]
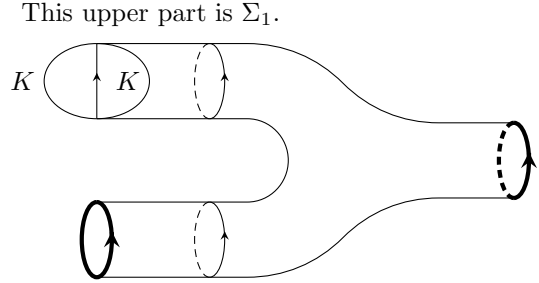
      {\textsc{Figure 11.} Open-closed cobordism homeomorphic to a closed window cobordism in figure 10.};
      
\end{tikzpicture}
\end{center}

Thus, for $a\in H_*(LM)$, 
\begin{equation*}
W_K(a)=\pm\mu_{\Sigma'}(a)=\pm\mu_{\Sigma_1}(1)\cdot a =\pm\chi(K)[c_0]\cdot a. 
\end{equation*}

This Proposition has the following immediate corollary. 

\begin{cor} [\textbf{Vanishing of double closed window operation}] \label{double closed window operation} For every oriented closed submanifolds $K,L$ of $M$, the composition of closed window operations $W_K$ and $W_L$ is identically zero. 
\begin{equation}
W_K\circ W_L\equiv 0: H_*(LM) \xrightarrow{W_L} H_*(LM) 
\xrightarrow{W_K}  H_*(LM).
\end{equation}
\end{cor}
\begin{proof}  For $a\in H_*(LM)$, we have 
\begin{equation*}
W_K\circ W_L(a)=W_K(\chi(L)[c_0]\cdot a)= \chi(K)\chi(L) [c_0]^2\cdot a=0, 
\end{equation*}
since $[c_0]^2=0$ in $H_*(LM)$ by a dimensional reason. 
\end{proof}

\begin{rem} So far we have discussed five types of string operations which identically vanish: open window operations, genus $1$ operations, handle attaching operations, double saddle operations, and double closed window operations. From open-closed cobordism point of view, these operations are closely related. Let $\Sigma$ be an open-closed cobordism for an open window operation. So $\Sigma$ is an annulus with one incoming open string and one outgoing open string along the outer boundary. Let $\Sigma_1$ and $\Sigma_2$ be two copies of $\Sigma$. We glue them in four different ways. See figure 12 for three of the four ways. 
\begin{enumerate}
\item[(i)] Glue $\Sigma_1$ and $\Sigma_2$ along outer free boundaries as well as along inner free boundary circles. We get an open-closed cobordism for genus $1$ operation. 
\item[(ii)] Glue $\Sigma_1$ and $\Sigma_2$ along inner boundary circles only, and label outer free boundaries by $I,J,K,L$. We get an open-closed cobordism for double saddle operation. 
\item[(iii)] Glue $\Sigma_1$ and $\Sigma_2$ along outer free boundaries, and label inner free boundary circles by $K,L$. We obtain an open-closed cobordism for double closed window operation. 
\item[(iv)] Glue $\Sigma_1$ and $\Sigma_2$ along inner free boundary circles and along one set of outer free boundaries, and label remaining free boundaries by $I,J$. We obtain an open-closed cobordism for handle attaching operation. 
\end{enumerate} 
All string operations associated to these open-closed cobordisms vanish. 
\end{rem}

\begin{figure}

\begin{center}
\begin{tikzpicture} [scale=0.7]
\path (0,2) coordinate (1);
\draw (1) -- ++(1,0) arc (270:315:1cm) arc (135:90:1cm) coordinate (2)
        -- ++(1,0) arc (90:45:1cm) arc (225:270:1cm) -- ++(1,0)
      -- ++(0.5,-1) -- ++(-1,0) arc (90:135:1cm) arc (315:270:1cm)                 coordinate (3) -- ++(-1,0) arc (270:225:1cm) arc (45:90:1cm) 
      -- ++(-1,0) -- cycle;
\draw[ultra thick] (1) -- ++(0.5,-1);
\path (1) ++(1,0) arc (270:315:1cm) arc (135:90:1cm)
       ++(1,0) arc (90:45:1cm) arc (225:270:1cm) ++(1,0) coordinate (8);
\draw[ultra thick] (8) -- ++(0.5,-1);
\path (2)++(0,-0.7) coordinate (4);
\path (1)++(4.2,-0.5) coordinate  (5);
\path (1)++(2.2,-0.5) coordinate (6);
\path (3)++(0,0.7) coordinate (7);
\draw (4)-- ++(1,0) to [out=0,in=120] (5) 
      to [out=300, in=0] (7) -- ++(-1,0) to [out=180,in=300]
      (6) to [out=120,in=180] (4);
 
\path (1)++(0,-1.8) coordinate(11);
\draw[ultra thick] (11) -- ++(0.5,-1);
\draw (11) -- ++(1,0) arc (270:315:1cm) arc (135:110:1cm);
\path (11) ++(1,0) arc (270:315:1cm) arc (135:110:1cm) coordinate (21);
\path (11) ++(1,0) arc (270:315:1cm) arc (135:90:1cm) coordinate (23);
\draw[densely dashed] (21) arc (110:90:1cm) -- ++(1,0) arc (90:45:1cm);
\path (21) arc (110:90:1cm) ++(1,0) arc (90:45:1cm) coordinate (22);
\draw (22) arc (225:270:1cm) -- ++(1,0)
      -- ++(0.5,-1) -- ++(-1,0) arc (90:135:1cm) arc (315:270:1cm)                 coordinate (31) -- ++(-1,0) arc (270:225:1cm) arc (45:90:1cm) 
      -- ++(-1,0) -- (11);
\path (22) arc (225:270:1cm) ++(1,0) coordinate (81);
\draw[ultra thick] (81) -- ++(0.5,-1);
\path (23)++(0,-0.7) coordinate (41);
\path (11)++(4.2,-0.5) coordinate  (51);
\path (11)++(2.2,-0.5) coordinate (61);
\path (31)++(0,0.7) coordinate (71);
\draw (41)-- ++(1,0) to [out=0,in=120] (51) 
      to [out=300, in=0] (71) -- ++(-1,0) to [out=180,in=300]
      (61) to [out=120,in=180] (41);

%Genus one cobordism. 
\path (13,8) coordinate (1a);
\draw (1a) arc (90:150:3cm and 2cm) 
      arc (330:270:1cm and 0.66cm) --++(-0.5,0);
\draw[densely dashed] (1a)++(0,-0.5) arc (90:150:3cm and 2cm) 
      arc (330:270:1cm and 0.66cm) --++(-0.3,0) coordinate (12a);
\draw (12a) -- ++(-0.4,0);
\draw (1a) arc (90:30:3cm and 2cm) arc (210:270:1cm and 0.66cm) 
      -- ++(0.5,0);
\draw[densely dashed] (1a)++(0,-0.5) arc (90:30:3cm and 2cm) arc (210:270:1cm and 0.66cm) 
      -- ++(0.27,0);
\draw (1a)++(0,-4) arc (270:210:3cm and 2cm) 
      arc (30:90:1cm and 0.66cm) -- ++(-0.5,0);
\draw (1a)++(0,-3.5) arc (270:210:3cm and 2cm) 
      arc (30:90:1cm and 0.66cm) -- ++(-0.27,0); 
\draw (1a)++(0,-4) arc (270:330:3cm and 2cm) arc (150:90:1cm and 0.66cm) 
      -- ++(0.5,0);
\draw (1a)++(0,-3.5) arc (270:330:3cm and 2cm) arc (150:90:1cm and 0.66cm) 
      -- ++(0.7,0);
\draw[densely dashed] (1a) arc (90:270:0.2cm and 0.7cm);
\path (1a) arc (90:270:0.2cm and 0.7cm) coordinate (2a);
\draw   (2a) arc (270:360:0.2cm and 0.7cm) 
        arc (0:90:0.2cm and 0.7cm);
\path (1a)++(0,-2.6) coordinate (3a);
\draw[densely dashed] (3a) arc (90:270:0.2cm and 0.7cm);
\draw    (1a)++(0,-4) arc (270:360:0.2cm and 0.7cm) 
arc (0:90:0.2cm and 0.7cm);
\draw[ultra thick] (1a)++(-3.95,-2) ellipse (0.25cm and 0.65cm);
\draw[ultra thick] (1a)++(3.95,-2.65) arc (270:360:0.25cm and 0.65cm)
          arc (0:90:0.25cm and 0.65cm);
\draw[densely dashed,ultra thick] 
       (1a)++(3.95,-1.35) arc (90:270:0.25cm and 0.65cm);
\path (1a)++(-1.5,-2) coordinate (4a);
\path (1a)++(1.5,-2) coordinate (5a);
\draw (4a) to [out=60,in=180] (2a) to [out=0,in=120] (5a);
\draw (4a)++(120:0.6cm) -- (4a) to [out=300,in=180] (3a) to [out=0,in=240] 
       (5a) -- ++(60:0.6cm);
\path (2a) ++ (0,0.5) coordinate  (6a);
\path (3a)++(0,-0.5) coordinate (7a);
\path (4a) ++(120:0.3cm) coordinate (8a);
\path (5a) ++(60:0.3cm) coordinate (9a);
\draw (8a) to [out=80, in=180] (6a) to [out=0,in=100] (9a);
\draw[densely dashed] (8a) to [out=260,in=180] (7a) to [out=0,in=280] (9a);

%Double Saddle cobordism
\path (11,2) coordinate (1b);
\draw (1b) -- ++(1,0) arc (270:315:1cm) arc (135:90:1cm) 
       -- ++(0.5,0) coordinate (2b)
        -- ++(0.5,0) arc (90:45:1cm) arc (225:270:1cm) 
        -- ++(1,0) coordinate (3b);
\draw[ultra thick] (3b) arc (180:195:0.3cm and 1.2cm);
\path (3b) arc (180:270:0.3cm and 1.2cm) coordinate (4b);
\path (3b) arc (180:345:0.3cm and 1.2cm) coordinate (5b);
\draw[ultra thick] (5b) arc (345:270:0.3cm and 1.2cm);
\draw[densely dashed,ultra thick] (4b) arc (270:195:0.3cm and 1.2cm);
\draw (5b) -- ++(-1,0) arc (90:135:1cm) arc (315:270:1cm) 
       -- ++(-0.5,0) coordinate (6b) -- ++(-0.5,0) arc (270:225:1cm) 
       arc (45:90:1cm) -- ++(-1,0) coordinate (7b);
\draw[ultra thick] (7b) arc (345:180:0.3cm and 1.2cm);
\path (7b) arc (345:270:0.3cm and 1.2cm) coordinate (8b);
\draw (8b) arc (90:0:1.5cm and 0.5cm) arc (360:270:1.5cm and 0.5cm);
\path (8b) arc (90:0:1.5cm and 0.5cm) coordinate (9b);
\draw (4b) arc (90:270:1.5cm and 0.5cm);
\path (4b) arc (90:180:1.5cm and 0.5cm) coordinate (10b);
\path (8b) arc (90:0:1.5cm and 0.5cm) 
      arc (360:270:1.5cm and 0.5cm) coordinate (19b);
\draw[ultra thick] (19b) arc (90:165:0.3cm and 1.2cm);
\path (19b) arc (90:165:0.3cm and 1.2cm) coordinate (11b);
\draw[ultra thick] (19b) arc (90:0:0.3cm and 1.2cm);
\draw[densely dashed] (11b)++(0.6,0) -- ++(0.4,0) arc (270:315:1cm) 
       arc (135:90:1cm) -- ++(0.5,0) coordinate (12b)
        -- ++(0.5,0) arc (90:45:1cm) arc (225:270:1cm) 
        -- ++(1,0) coordinate (13b);
\draw (11b) -- ++(0.6,0);
\draw[densely dashed,ultra thick] (13b) arc (165:90:0.3cm and 1.2cm);
\path (13b) arc (165:0:0.3cm and 1.2cm) coordinate (14b);
\draw[ultra thick] (14b) arc (0:90:0.3cm and 1.2cm);
\draw (14b) -- ++(-1,0) arc (90:135:1cm) arc (315:270:1cm) 
       -- ++(-0.5,0) coordinate (15b) -- ++(-0.5,0) arc (270:225:1cm) 
       arc (45:90:1cm) -- ++(-1,0) coordinate (16b);
\draw (9b) arc (180:360:1.41cm and 0.4cm);
\draw[densely dashed] (9b) arc (180:0:1.41cm and 0.4cm);
\path (9b) arc (180:275:1.41cm and 0.4cm) coordinate (18b);
\path (9b) arc (180:95:1.41cm and 0.4cm) coordinate (17b);
\draw (6b) to [out=250,in=90] (18b);
\draw (18b) to [out=270,in=110] (15b);
\path (17b)++(0,0.5) coordinate (21b);
\draw (2b) to [out=290,in=90] (21b);
\draw[densely dashed] (21b) -- (17b);
\draw[densely dashed] (17b) to [out=270,in=70] (12b);

%Double closed window cobordism. 
\path (13,-4) coordinate (1c);
\draw (1c) arc (90:150:3cm and 2cm) 
      arc (330:270:1cm and 0.66cm) --++(-0.5,0);
\draw[densely dashed] (1c)++(0,-0.5) arc (90:150:3cm and 2cm) 
      arc (330:270:1cm and 0.66cm) --++(-0.3,0) coordinate (12c);
\draw (12c) -- ++(-0.4,0);
\draw (1c) arc (90:30:3cm and 2cm) arc (210:270:1cm and 0.66cm) 
      -- ++(0.5,0);
\draw[densely dashed] (1c)++(0,-0.5) arc (90:30:3cm and 2cm) arc (210:270:1cm and 0.66cm) 
      -- ++(0.27,0);
\draw (1c)++(0,-4) arc (270:210:3cm and 2cm) 
      arc (30:90:1cm and 0.66cm) -- ++(-0.5,0);
\draw (1c)++(0,-3.5) arc (270:210:3cm and 2cm) 
      arc (30:90:1cm and 0.66cm) -- ++(-0.27,0); 
\draw (1c)++(0,-4) arc (270:330:3cm and 2cm) 
      arc (150:90:1cm and 0.66cm) 
      -- ++(0.5,0);
\draw (1c)++(0,-3.5) arc (270:330:3cm and 2cm) 
       arc (150:90:1cm and 0.66cm) 
      -- ++(0.7,0);
\draw[ultra thick] (1c)++(-3.95,-2) ellipse (0.25cm and 0.65cm);
\draw[ultra thick] (1c)++(3.95,-2.65) 
    arc (270:360:0.25cm and 0.65cm) arc (0:90:0.25cm and 0.65cm);
\draw[densely dashed,ultra thick] 
       (1c)++(3.95,-1.35) arc (90:270:0.25cm and 0.65cm);
\draw (1c)++(0,-1.6) ellipse (1.8cm and 0.7cm);
\draw[densely dashed] (1c)++(0,-2.4) ellipse (1.8cm and 0.7cm);

\path (1c) arc (90:95:3cm and 2cm) coordinate (2c);
\path (1c)++(0,-4) arc (270:275:3cm and 2cm) coordinate (3c);

\draw[densely dashed] (2c) arc (90:210:0.4cm and 0.7cm) 
     -- ++(285:0.68cm);
\draw  (2c) arc (90:30:0.4cm and 0.7cm) -- ++(285:0.58cm);
     
\draw (3c) arc (270:360:0.4cm and 0.7cm) arc (0:30:0.4cm and 0.7cm) 
       -- ++(105:0.68cm);
\draw[densely dashed] (3c) arc (270:210:0.4cm and 0.7cm) 
      -- ++(105:0.58cm);

\draw[->,>=stealth,thick] (6,3) -- node[above left,text width=4.5cm] 
{(i) Glue along both outer and inner free boundaries.} (9,5) ;
\draw[->,>=stealth,thick] (7,0.5) -- node[above, text width=2.5cm] {(ii) Glue along the inner free boundary.}(10.5,0.5) ;
\draw[->,>=stealth,thick] (6,-2) -- node[below left, text width=4cm] {(iii) Glue along the outer free boundaries.} (9,-5) ;

\draw (14,3.7) node[below, text width=7cm] {(i) Genus $1$ open-closed cobordism.};
\draw (14,-2.5) node[below,text width=7cm] {(ii) Double saddle open-closed cobordism.};
\draw (13,-8.5) node[below,text width=6cm] {(iii) Double closed window open-closed cobordism.};

\draw (0,-10) node[below right, text width=14cm] {\textsc{Figure 12.} Three of four ways to glue two copies of open window cobordisms giving rise to three basic open-closed cobordisms with vanishing string operations.};

\end{tikzpicture}
\end{center}

\end{figure}

\bigskip

\section{Proof of Theorems} 
%Section 4

In this section, we will combine vanishing results in the previous section with some auxiliary results to prove the main theorem. Let $\Sigma$ be a connected open-closed cobordism, and let $\mu_{\Sigma}$ be the associated orientable string operation determined up to a sign. 

The first proposition proves the first part of the main theorem. The meaning of integer invariants $g(\Sigma), \omega(\Sigma), p(\Sigma), q(\Sigma), r(\Sigma), s(\Sigma), t(\Sigma)$ associated to the open-closed cobordism $\Sigma$ is described in the introduction. 

\begin{prop}  Let $\Sigma$ be a connected open-closed cobordism. The string operation $\mu_{\Sigma}$ identically vanishes if $\Sigma$ satisfies one of the following conditions. 
\begin{gather*}
\textup{(i)} \ \ g(\Sigma)\ge1. \qquad 
\textup{(ii)} \ \ \omega(\Sigma)\ge2. \qquad 
\textup{(iii)} \ \ t(\Sigma)\ge2. \\
\textup{(iv)} \ \ q(\Sigma)\ge 3. \qquad 
\textup{(v)} \ \ s(\Sigma)\ge1 \text{ and } s(\Sigma)+q(\Sigma)\ge2. 
\end{gather*}
In \textup{(v)}, assume that at least one label $K$ of outgoing open strings has dimension $\dim K<d$. 
\end{prop} 
\begin{proof} The conditions (i) and (ii) come from Proposition \ref{open window operation}, Proposition \ref{genus one operation}, Proposition \ref{handle attaching operation}, and Corollary \ref{double closed window operation}. Proposition \ref{handle attaching operation} is needed for (i) to cover the case when $\Sigma$ has a boundary containing both incoming and outgoing open strings ($t(\Sigma)\ge1$). When a connected open-closed cobordism $\Sigma$ has at least two boundaries containing both incoming and outgoing open strings ($t(\Sigma)\ge2$), $\Sigma$ contains an open-closed cobordism surface for double saddle operation. Thus, the condition (iii) comes from Proposition \ref{double saddle operation}. In \cite{T3}, we showed that if $\Sigma$ is a connected closed string cobordism with at least three outgoing closed strings, then the associated string operation identically vanishes. The condition (iv) comes from this fact. 

For (v), let $\Sigma'$ be a genus $0$ surface with three boundary circles having one incoming closed string, one outgoing closed string, and one outgoing open string whose end points are labeled by the closed oriented submanifold $K$ of dimension $\dim K=k<d$. When $\Sigma$ satisfies the condition (v), it contain an open-closed cobordism homeomorphic to $\Sigma'$. The string operation associated to $\Sigma'$ is 
\begin{equation*}
\mu_{\Sigma'}=(\theta_K\otimes 1)\varphi: H_*(LM) \longrightarrow H_*(LM)\otimes H_*(P_{KK}), 
\end{equation*}
where for $a\in H_*(LM)$, we have $(\theta_K\otimes 1)\varphi(a)=(\theta_K\otimes 1)\bigl([c_0]\otimes([c_0]\cdot a)\bigr)$. Since the operation $\theta_K$ lowers degree by $d-k>0$, $\theta([c_0])=0\in H_*(P_{KK})$. Hence $\mu_{\Sigma'}\equiv 0$. Thus, string operation associated to connected open-closed cobordisms satisfying the condition (v) identically vanish. 
\end{proof}

We use genus $g(\Sigma)$ and window numbers $\omega(\Sigma)$ as primary invariants to classify open-closed cobordisms with vanishing string operations. 

\begin{prop} Let $\Sigma$ be a connected open-closed cobordism with the genus $g(\Sigma)=0$ and the window number $\omega(\Sigma)=1$. 
\begin{enumerate}
\item If $\Sigma$ has an incoming or an outgoing open string, then $\mu_{\Sigma}\equiv 0$. 
\item If $\Sigma$ has no open strings but has at least $2$ outgoing closed strings, then $\mu_{\Sigma}\equiv0$. 
\end{enumerate}
\end{prop}
\begin{proof} (1) If $\Sigma$ has either an incoming or outgoing open string ($r+s+t\ge1$), then $\Sigma$ contains an open-closed cobordism for the open window operation. Hence the vanishing of the string operation $\mu_{\Sigma}$ is a consequence of Proposition \ref{open window operation}. Thus, for a surface $\Sigma$ with $g=0$ and $\omega=1$, nontriviality of the associated string operation requires that $\Sigma$ has no open strings. 

For (2), let $\Sigma'$ be an open-closed cobordism with one incoming closed string, two outgoing closed strings and one free boundary circle labeled by a closed oriented submanifold $K$. Then the associated string operation $\mu_{\Sigma'}$ lowering degree by $2d$ is given by 
\begin{equation*}
\mu_{\Sigma'}=\varphi\circ W_K: H_*(LM) \longrightarrow H_*(LM)\otimes H_*(LM).
\end{equation*}
For $a\in H_*(LM)$, we have $\mu_{\Sigma'}(a)=\varphi(\chi(K)[c_0]\cdot a)=\chi(M)\chi(K)[c_0]\otimes ([c_0]^2\cdot a)=0$, since $[c_0]^2=0$ in $H_*(LM)$. Thus the string operation $\mu_{\Sigma'}$ identically vanish. Hence the original string operation $\mu_{\Sigma}$ must also vanish. This completes the proof. 
\end{proof}

This proves the part (II) of Theorem A. Since in string topology, there must be at least one outgoing open or closed string, part (2) of the above proposition leaves the case of closed cobordism with one outgoing closed string and one closed window. The associated string operation can be nontrivial. This is the first type of open-closed cobordisms with possibly nontrivial string operations listed in Theorem B.  

Next we consider the case in which $g(\Sigma)=0$ and $\omega(\Sigma)=0$, and we prove part (III) of Theorem A. Recall that the number of boundary components $t(\Sigma)$ containing both incoming and outgoing open strings must be at most one for nontrivial string operations by the first part of Theorem A.  

\begin{prop} Let $\Sigma$ be a connected genus $0$ open-closed cobordism with no windows. Suppose $\Sigma$ has exactly one boundary containing both incoming and outgoing open strings \textup{(}that is, $t(\Sigma)=1$\textup{)}, with a label of dimension less than $d=\dim M$. If $\Sigma$ contains an outgoing closed string or a boundary component containing only outgoing open strings, then the associated string operation identically vanish. Namely, if $q(\Sigma)+s(\Sigma)\ge1$, then $\mu_{\Sigma}\equiv 0$. 
\end{prop}
\begin{proof} Such an open-closed cobordism $\Sigma$ contains a connected genus $0$ open-closed cobordism $\Sigma'$ with one outgoing closed string and one boundary containing one incoming open string and one outgoing open string. Let the labels of the free boundary components be $I,J$ such that $\dim I<\dim M$. See figure 13. The associated string operation is of the form
\begin{equation*}
\mu_{\Sigma'}=\varphi_{IJ}: H_*(P_{IJ}) \longrightarrow H_*(LM)\otimes H_*(P_{IJ}).
\end{equation*}
where the $H_*(LM)$ comodule map $\varphi_{IJ}=j_*\iota_!$ lowering degree by $d$ is defined using maps $j$ and $\iota$ given in the following diagram:
\begin{equation*}
\begin{CD} 
P_{IJ} @<{\iota}<< LM\underset{M}{\times}P_{IJ} @>{j}>> LM\times P_{IJ}  \\
@V{(p_{\frac14},p_{\frac34})}VV  @V{p}VV  @V{(p\times p_{\frac12})}VV \\
M\times M @<{\phi}<< M  @>{\phi}>> M\times M ,
\end{CD}
\end{equation*} 
where $LM\times_MP_{IJ}$ consists of pairs $(\gamma,\eta)\in LM\times P_{IJ}$ such that $\gamma(0)=\eta(\frac12)$. We show that the comodule map $\varphi_{IJ}$ is identically zero.

\begin{figure}

\begin{center}
\begin{tikzpicture}
\draw (1,0) arc (180:90:2cm) -- ++(1,0) arc (90:270:0.2cm and 0.5cm) arc (270:360:0.2cm and 0.5cm) arc (0:90:0.2cm and 0.5cm) ++(0,-1) -- ++(-1,0) arc (90:180:1cm) arc (360:180:0.5cm and 0.2cm);
\path (1,0) arc (180:90:2cm) -- ++(1,0) coordinate (1);
\draw[ultra thick] (1) arc (90:360:0.2cm and 0.5cm) arc (0:90:0.2cm and 0.5cm);
\draw[->,>=stealth,ultra thick] (1)++(-0.2,-0.5) -- ++(0,0.02);
\draw[densely dashed] (1,0) arc (180:0:0.5cm and 0.2cm);
\draw (1.1,0.7) -- (0.3,0.7) --(-0.3,-0.7) -- (3.7,-0.7) -- (4.3,0.7) -- (2.3,0.7);
\draw[ultra thick] (0.3,0.7) --(-0.3,-0.7) (3.7,-0.7) -- (4.3,0.7);
\draw[->,>=stealth,ultra thick] (0,0) -- ++(0.03,0.07);
\draw[->,>=stealth,ultra thick]  (4,0) -- ++(0.03,0.07);
\draw (0.3,0.7) node[above right] {$J$};
\draw (-0.3,-0.7) node[below right] {$I$};
\draw (3.7,-0.7) node[below left] {$I$};
\draw (4.3,0.7) node[above] {$J$};
\draw (2,-1.2) node[below, text width=10cm] 
 {\textsc{Figure 13}. Open-closed cobordism for $H_*(LM)$-comodule structure in $H_*(P_{IJ})$ : 
 $H_*(P_{IJ}) \longrightarrow H_*(LM)\otimes H_*(P_{IJ})$};
\end{tikzpicture}
\end{center}

\end{figure}

We can easily verify that the comodule map $\varphi_{IJ}$ satisfies the property 
\begin{equation*}
\varphi_{IJ}(a\cdot_I b)=\varphi_{II}(a)\cdot_I b
\end{equation*}
for any $a\in H_*(P_{II})$ and $b\in H_*(P_{IJ})$. Letting $a=[I]$, the unit in $H_*(P_{II})$, by degree reason we have $\varphi_{II}([I])=0$, since $\varphi_{II}$ lowers the degree by $d$ and by assumption $\dim I<d$. Hence for any $b\in H_*(P_{IJ})$, $\varphi_{IJ}(b)=\varphi_{IJ}([I]\cdot_Ib)
=\varphi_{II}([I])\cdot_I b=0$. Hence the string operation $\mu_{\Sigma'}=\varphi_{IJ}$ identically vanishes. This implies that the original string operation $\mu_{\Sigma}$ must vanish, too. 
\end{proof} 

The above proposition leaves the case in which the genus $0$ connected open-closed cobordisms $\Sigma$ with no window have no closed outgoing strings nor boundaries with only outgoing open strings ($q(\Sigma)+s(\Sigma)=0$). The associated string operations can be nontrivial. This cobordism is listed in the introduction as the fifth case in Theorem B.  The three cases (2), (3) and (4) in Theorem B correspond to those open-closed cobordisms excluded in part (IV) of Theorem A. This completes the proof of Theorem B.

\bibliography{bibliography}
\bibliographystyle{plain}

\end{document}